# HAUSDORFF MEASURE OF ARCS AND BROWNIAN MOTION ON BROWNIAN SPATIAL TREES

By David A. Croydon

*University of Warwick*

A Brownian spatial tree is defined to be a pair $(\mathcal{T}, \phi)$, where $\mathcal{T}$ is the rooted real tree naturally associated with a Brownian excursion and $\phi$ is a random continuous function from $\mathcal{T}$ into $\mathbb{R}^d$ such that, conditional on $\mathcal{T}$, $\phi$ maps each arc of $\mathcal{T}$ to the image of a Brownian motion path in $\mathbb{R}^d$ run for a time equal to the arc length. It is shown that, in high dimensions, the Hausdorff measure of arcs can be used to define an intrinsic metric $d_\mathcal{S}$ on the set $\mathcal{S} := \phi(\mathcal{T})$. Applications of this result include the recovery of the spatial tree $(\mathcal{T}, \phi)$ from the set $\mathcal{S}$ alone, which implies in turn that a Dawson–Watanabe super-process can be recovered from its range. Furthermore, $d_\mathcal{S}$ can be used to construct a Brownian motion on $\mathcal{S}$, which is proved to be the scaling limit of simple random walks on related discrete structures. In particular, a limiting result for the simple random walk on the branching random walk is obtained.

**1. Introduction.** Super-processes are measure-valued diffusions that arise naturally as the scaling limits of discrete branching particle models in Euclidean space; see [26] for an introduction to this area. Describing the genealogy of super-processes provided one of the original motivations for the study of continuous branching structures, which has been intense in recent years; [22] is an up-to-date survey article. The second key component in defining a super-process is the description of how "particles" proceed through $\mathbb{R}^d$. A particularly important example of a super-process is the Dawson–Watanabe super-process, $(Y_t)_{t \geq 0}$, say which has a binary branching structure and whose spatial motion is given by Brownian motion in $\mathbb{R}^d$ (see Section 3 for a precise definition). If $d \geq 2$, it is known that, for each fixed $t > 0$, the measure $Y_t$ can almost surely be obtained from its support $S(Y_t)$ as a Hausdorff measure on this set [12, 23]; when $d \geq 3$, another representation of $Y_t$ is provided in [28]. Similarly, when $d \geq 4$, for almost-every realization of the super-process, it is possible to reconstruct the total occupation measure $\int_0^\infty Y_t \, dt$ of the super-process from its range $\mathcal{R}_1$, which is defined below at (3.2), as a Hausdorff measure [20]. In high dimensions, by defining an intrinsic metric on the range of the super-process, we will show that it is possible to separate the branching structure and spatial motion of the super-process and, thereby, extend these results. Specifically, we are able to prove that, for $d \geq 8$, the Dawson–Watanabe super-process $(Y_t)_{t \geq 0}$ can almost surely be reconstructed from knowledge of its range $\mathcal{R}_1$ alone (see Corollary 5.4).

The framework for this article is the space of spatial trees introduced by Duquesne and Le Gall in [14]. In particular, we will consider Brownian spatial trees, by which we mean pairs of the form $(\mathcal{T}, \phi)$, where $\mathcal{T}$ is the rooted real tree naturally associated with a Brownian excursion and $\phi$ is a random continuous function from $\mathcal{T}$ into $\mathbb{R}^d$ such that, conditional on $\mathcal{T}$, $\phi$ maps each arc of $\mathcal{T}$ to the image of a Brownian motion path in $\mathbb{R}^d$ run for a time equal to the arc length (see Sections 2.1 and 2.2 for details). The key step in establishing the result described at the end of the previous paragraph is showing that, when $d \geq 8$, the set $\mathcal{S} := \phi(\mathcal{T})$ almost surely determines the spatial tree $(\mathcal{T}, \phi)$ (see Corollary 5.3), and to do this, we apply two main ideas. First, we use known intersection properties of super-processes to check that, when $d \geq 8$, the map $\phi : \mathcal{T} \to \mathcal{S}$ is a homeomorphism (see Section 3). It follows that $\mathcal{S}$ is almost surely a dendrite (an arcwise-connected topological space containing no subset homeomorphic to the circle) and, therefore, between any two points of $\mathcal{S}$ there is a unique arc in $\mathcal{S}$. Second, it was proved in [7] that, when $d \geq 3$, a Brownian motion path in $\mathbb{R}^d$ run for a time $t$ has Hausdorff measure $t$, almost surely, with respect to the measure function $c_d x^2 \ln \ln x^{-1}$, where $c_d$ is a deterministic constant that depends only upon $d$. Since arcs in $\mathcal{S}$ are, by construction, segments of Brownian motion paths, we can combine these two observations to define a metric $d_\mathcal{S}$ on $\mathcal{S}$ by setting, for $x_1, x_2 \in \mathcal{S}$, the distance $d_\mathcal{S}(x_1, x_2)$ to be equal to the Hausdorff measure, with respect to the measure function $c_d x^2 \ln \ln x^{-1}$, of the unique arc between $x_1$ and $x_2$ in $\mathcal{S}$. It is then possible to demonstrate that $\phi$ is actually an isometry from $(\mathcal{T}, d_\mathcal{T})$, where $d_\mathcal{T}$ is the natural metric on $\mathcal{T}$, to $(\mathcal{S}, d_\mathcal{S})$ almost surely, and, consequently, we obtain that $(\mathcal{T}, \phi)$ and $(\mathcal{S}, I)$ are equivalent spatial trees almost surely, where $I$ is the identity map on $\mathbb{R}^d$ (see Proposition 5.2).

A second application of the metric $d_\mathcal{S}$ is that it allows the construction of a natural diffusion on the set $\mathcal{S}$. First note that, by applying results of [18], a Dirichlet form can be constructed on any compact real tree equipped with a suitable Borel measure, and following the arguments of [9], Section 8, it is possible to check that the corresponding diffusion is actually, a Brownian motion on the relevant space, as defined by Aldous in [2]. Since $(\mathcal{S}, d_\mathcal{S})$ is a real tree, it fits naturally into this setting and, therefore, to define a Brownian motion $X^{\mathcal{S}} = (X_t^{\mathcal{S}})_{t\geq 0}$ on $\mathcal{S}$, it remains to choose an appropriate measure. For $d \geq 8$, the canonical measure on $\mathcal{S}$, which we will denote by $\mu^{\mathcal{S}}$, is equal to the Hausdorff measure on $\mathcal{S}$ with measure function proportional to $x^4 \ln \ln x^{-1}$ ([20], Theorem 6.1), and can also be interpreted as $\mu^{\mathcal{T}} \circ \phi^{-1}$, where $\mu^{\mathcal{T}}$ is the natural measure on $\mathcal{T}$ (see Sections 2.1 and 2.2). Consequently, using the fact that $\phi$ is an isometry, it is possible to show that the resulting process $X^{\mathcal{S}}$ can also be written as $\phi(X^{\mathcal{T}})$ almost surely, where $X^{\mathcal{T}}$ is the Brownian motion associated with $\mathcal{T}$ and $\mu^{\mathcal{T}}$ (see Proposition 6.1); we observe that $\phi(X^{\mathcal{T}})$ is actually defined for any $d \geq 1$. In addition to defining the laws of $\phi(X^{\mathcal{T}})$ and $X^{\mathcal{S}}$ for almost-every realization of $(\mathcal{T}, \phi)$ and $\mathcal{S}$ respectively, which are the quenched versions of the laws, by adapting the arguments of [8], we demonstrate the measurability of the construction, which allows us to define related annealed laws (where we average out over all realizations of the spatial trees).

While we do not pursue it in depth here, let us remark that the representation $\phi(X^{\mathcal{T}})$ of $X^{\mathcal{S}}$, where $\phi$ is an isometry, means that we are immediately able to deduce many properties of the Brownian motion on $\mathcal{S}$ from known results about $X^{\mathcal{T}}$. For example, it follows from results appearing in [9] that, when $d \geq 8$, the diffusion $X^{\mathcal{S}}$ on almost-every realization of $\mathcal{S}$ admits a transition density $(p_t^{\mathcal{S}}(x,y))_{x,y\in\mathcal{S},t>0}$ that satisfies

$$\lim_{t\to 0} \frac{2 \ln p_t^{\mathcal{S}}(x,x)}{-\ln t} = \frac{4}{3} \qquad \forall x \in \mathcal{S}.$$

Using the terminology of the diffusion on fractal literature, this result could be interpreted as a version of the statement that the spectral dimension of the Brownian motion on $\mathcal{S}$ is $4/3$, almost surely. More detailed transition density asymptotics are obtained in [9]. Furthermore, asymptotics for the spectrum of the generator of the diffusion $X^{\mathcal{S}}$ are consequences of the results appearing in [10].

One reason for wanting to define a canonical process on the set $\mathcal{S}$ is that it provides an archetype for the scaling limit of the simple random walks on the graph-based models that converge in some sense to the integrated super-Brownian excursion, which was originally defined in [4] and can be thought of as the measure $\mu^{\mathcal{S}}$ conditioned to have total mass equal to one. Examples of discrete models which fall into this category include conditioned branching random walks, lattice trees in high dimensions and large critical percolation clusters in high dimensions (see [27], Section 16, for background). To prove a first result in this direction, we consider a sequence $\{(T_n, \phi_n)\}_{n\geq 1}$ of random "graph spatial trees," by which we mean that $T_n$ is a random (rooted) ordered graph tree and $\phi_n$ is a random embedding of $T_n$ into $\mathbb{R}^d$. Our main assumption is that, for each $n$, $T_n$ has $n$ vertices, and also that the discrete tours associated with $(n^{-1/2}T_n, n^{-1/4}\phi_n)$ converge to the normalized Brownian tour; in both the discrete and continuous cases, a tour is a continuous function that encapsulates the branching and spatial motion of the relevant spatial trees (see Sections 2 and 8 for exact definitions). Under the appropriate versions of this condition, we are able to deduce quenched and annealed versions of the statement that the process

$$(n^{-1/4}\phi_n(X^{T_n}_{n^{3/2}t}))_{t\geq 0},$$

where $X^{T_n}$ is the usual discrete time simple random walk on the vertices of $T_n$ started from the root, converges to the process $\phi(X^{\mathcal{T}})$, which, as remarked above, is identical to the process $X^{\mathcal{S}}$ in high dimensions (see Theorems 8.1 and 9.1). To prove these results, we apply ideas from [8], which demonstrates the convergence of $(n^{-1/2}X^{T_n}_{n^{3/2}t})_{t\geq 0}$ to $X^{\mathcal{T}}$ under a related assumption that does not include any spatial component. Branching random walks with a critical offspring distribution that decays exponentially at infinity and a step distribution that satisfies an $o(x^{-4})$ tail bound, conditioned on the total number of offspring, are a special case of graph spatial trees known to satisfy assumptions that allow the above scaling limit results to be applied (see Section 10 for details).

The article is arranged as follows. In Section 2 we introduce much of the notation for real trees, spatial trees and tours that will be used throughout the article. Section 3 contains a proof of the fact that $\phi:\mathcal{T}\to\mathcal{S}$ is a homeomorphism in high dimensions, and in Section 4 we investigate the Hausdorff measure of arcs of $\mathcal{S}$. The first half of the article is concluded in Section 5, where we define $d_{\mathcal{S}}$, verify that $\phi:(\mathcal{T},d_{\mathcal{T}})\to(\mathcal{S},d_{\mathcal{S}})$ is an isometry in high dimensions and prove the super-process result described in our opening paragraph. The second half of the article is devoted to the study of $X^{\mathcal{S}}$ and $\phi(X^{\mathcal{T}})$. We first define the quenched and annealed laws of these processes in Sections 6 and 7, respectively. The quenched and annealed convergence results for simple random walks on graph spatial trees are then proved in Sections 8 and 9, respectively. Finally, in Section 10 we apply these results to the simple random walk on the branching random walk.

## 2. Notation and preliminaries.

2.1. *Real trees and excursions.*  At the core of our discussion will be the collection of metric spaces known as real trees, for which we use the following definition. Note that much of the notation and many of the definitions used in this section are borrowed from [14] and other works by the same authors.

DEFINITION 2.1.  A metric space $(\mathcal{T},d_{\mathcal{T}})$ is a real tree if the following properties hold for every $\sigma_1,\sigma_2\in\mathcal{T}$:

(a) There is a unique isometric map $\gamma^{\mathcal{T}}_{\sigma_1,\sigma_2}$ from $[0,d_{\mathcal{T}}(\sigma_1,\sigma_2)]$ into $\mathcal{T}$ such that $\gamma^{\mathcal{T}}_{\sigma_1,\sigma_2}(0)=\sigma_1$ and $\gamma^{\mathcal{T}}_{\sigma_1,\sigma_2}(d_{\mathcal{T}}(\sigma_1,\sigma_2))=\sigma_2$.

(b) If $\gamma$ is a continuous injective map from $[0,1]$ into $\mathcal{T}$ such that $\gamma(0) = \sigma_1$ and $\gamma(1) = \sigma_2$, then $\gamma([0,1]) = \gamma^{\mathcal{T}}_{\sigma_1,\sigma_2}([0, d_{\mathcal{T}}(\sigma_1, \sigma_2)])$.

A rooted real tree is a real tree $(\mathcal{T}, d_{\mathcal{T}})$ with a distinguished vertex $\rho = \rho(\mathcal{T})$ called the root.

All the real trees we consider will be rooted, although for brevity we will often write simply $\mathcal{T}$ to represent $(\mathcal{T}, d_{\mathcal{T}}, \rho)$. The arc between two vertices $\sigma_1$ and $\sigma_2$ of a real tree $\mathcal{T}$ will be denoted by $\Gamma^{\mathcal{T}}_{\sigma_1,\sigma_2}$; more specifically, $\Gamma^{\mathcal{T}}_{\sigma_1,\sigma_2}$ is the image of $\gamma^{\mathcal{T}}_{\sigma_1,\sigma_2}$. An observation that will be useful to us is that between any three points $\sigma_1, \sigma_2, \sigma_3$ of a real tree $\mathcal{T}$ there is a unique branch-point $b^{\mathcal{T}}(\sigma_1, \sigma_2, \sigma_3) \in \mathcal{T}$ that satisfies

$$\text{(2.1)} \qquad \{b^{\mathcal{T}}(\sigma_1, \sigma_2, \sigma_3)\} = \Gamma^{\mathcal{T}}_{\sigma_1,\sigma_2} \cap \Gamma^{\mathcal{T}}_{\sigma_2,\sigma_3} \cap \Gamma^{\mathcal{T}}_{\sigma_3,\sigma_1}.$$

A useful decomposition of a real tree is provided by the subsets containing points equidistant from the root. In particular, define the subset of $\mathcal{T}$ at level $t$ to be

$$\text{(2.2)} \qquad \mathcal{T}_t := \{\sigma \in \mathcal{T} : d_{\mathcal{T}}(\rho, \sigma) = t\}.$$

The height of a real tree is given by $h(\mathcal{T}) := \sup\{d_{\mathcal{T}}(\rho, \sigma) : \sigma \in \mathcal{T}\}$, and we clearly have $\mathcal{T}_t = \varnothing$ for $t > h(\mathcal{T})$. We will also be interested in the decomposition of a real tree $\mathcal{T}$ into the subset below a certain level and the collection of subtrees of $\mathcal{T}$ that start at this level. To introduce this, define the set

$$\text{(2.3)} \qquad \text{tr}_t(\mathcal{T}) := \{\sigma \in \mathcal{T} : d_{\mathcal{T}}(\rho, \sigma) \leq t\}$$

which is the truncation of the real tree $\mathcal{T}$ at level $t$. Furthermore, let $\mathcal{T}^{i,o}$, $i \in \mathcal{I}_t$, be the connected components of the open set $\mathcal{T} \setminus \text{tr}_t(\mathcal{T})$. Note that if $h(\mathcal{T}) \leq t$, the collection $\mathcal{I}_t$ is empty. Observe that the ancestor of $\sigma$ at level $t$ [the unique point on the arc between $\rho$ and $\sigma$ with $d_{\mathcal{T}}(\rho, \sigma) = t$] must be the same for each $\sigma \in \mathcal{T}^{i,o}$, and we will denote it by $\rho^i$. Now define $\mathcal{T}^i := \mathcal{T}^{i,o} \cup \{\rho^i\}$, which is a real tree when endowed with the metric induced by $d_{\mathcal{T}}$ and we set its root to be $\rho^i$.

Of course, there are collections of real trees that are indistinguishable as metric spaces. We will denote by $\mathbb{T}$ the set of equivalence classes of compact rooted real trees, under the assumption that two rooted real trees are equivalent if and only if there exists a root preserving isometry between them. The set $\mathbb{T}$ can be equipped with the (pointed) Gromov–Hausdorff distance, $d_{\text{GH}}$, say, and it has been proved that $(\mathbb{T}, d_{\text{GH}})$ is a Polish space; see [16], Theorem 1. In our discussion of the properties of elements of $\mathbb{T}$ it will suffice to consider one particular real tree of each equivalence class. For detailed remarks about the technicalities of defining objects such as local times as we do below, see [14].

A particularly useful representation of the real trees that are studied in this article is provided through excursions. Define the space of excursions, $\mathcal{V}$, to be the set of continuous functions $v: \mathbb{R}_+ \to \mathbb{R}_+$ for which there exists a $\tau(v) \in (0, \infty)$ such that $v(t) > 0$ if and only if $t \in (0, \tau(v))$. Given a function $v \in \mathcal{V}$, we define a distance on $[0, \tau(v)]$ by setting

$$(2.4) \qquad d_v(s,t) := v(s) + v(t) - 2m_v(s,t),$$

where $m_v(s,t) := \inf\{v(r) : r \in [s \wedge t, s \vee t]\}$, and then use the equivalence

$$(2.5) \qquad s \sim t \quad \Leftrightarrow \quad d_v(s,t) = 0,$$

to define $\mathcal{T}_v := [0, \tau(v)]/\sim$. Denoting by $[s]$ the equivalence class containing $s$, it is elementary (see [14], Section 2) to check that $d_{\mathcal{T}_v}([s],[t]) := d_v(s,t)$ defines a metric on $\mathcal{T}_v$, and also that with this metric $\mathcal{T}_v$ is a real tree. The root of the tree $\mathcal{T}_v$ is defined to be the equivalence class $[0]$, and is denoted by $\rho_v$. A natural volume measure to impose upon $\mathcal{T}_v$ is the projection of the Lebesgue measure on $[0, \tau(v)]$. In particular, for open $A \subseteq \mathcal{T}_v$, let

$$(2.6) \qquad \mu_v(A) := \lambda(\{t \in [0, \tau(v)] : [t] \in A\}),$$

where $\lambda$ is the usual one-dimensional Lebesgue measure. This defines a Borel measure on $(\mathcal{T}_v, d_{\mathcal{T}_v})$, with total mass equal to $\tau(v)$. We will usually suppress the dependence on $v$ from the notation for all of these objects when it is clear which excursion is being considered.

We will be interested in the measure $\Theta$ on $\mathbb{T}$, which is defined by

$$(2.7) \qquad \Theta(A) := N(\{v : \mathcal{T}_v \in A\})$$

for measurable $A \subseteq \mathbb{T}$, where $N$ is the usual Itô excursion measure, normalized so that the tail of the height of a tree chosen from $\Theta$ is given by

$$(2.8) \qquad \Theta(h(\mathcal{T}) > \varepsilon) = \varepsilon^{-1}.$$

A real tree $\mathcal{T}$ chosen according to $\Theta$ is an example of a Lévy tree, as introduced in [14], and an important result of [14] is that a Lévy tree admits an intrinsic "local time" measure, as described by the following theorem.

THEOREM 2.2 ([14], Theorem 4.2). *For every $t \geq 0$ and $\Theta$-a.e. $\mathcal{T} \in \mathbb{T}$, we can define a finite measure $\ell_t$ on $\mathcal{T}$ in such a way that the following properties hold:*

(a) *$\ell_0 = 0$ and, for every $t > 0$, $\ell_t$ is supported on $\mathcal{T}_t$, $\Theta$-a.e.*
(b) *for every $t > 0$, $\{\ell_t \neq 0\} = \{h(\mathcal{T}) > t\}$, $\Theta$-a.e.*
(c) *for every $t > 0$, we have $\Theta$-a.e. for every bounded continuous function $\varphi$ on $\mathcal{T}$,*

$$\int_{\mathcal{T}} \varphi \, d\ell_t = \lim_{\varepsilon \to 0} \sum_{i \in \mathcal{I}_t} \varepsilon \varphi(\rho^i) \mathbf{1}_{\{h(\mathcal{T}^i) \geq \varepsilon\}} = \lim_{\varepsilon \to 0} \sum_{i \in \mathcal{I}_{t-\varepsilon}} \varepsilon \varphi(\rho^i) \mathbf{1}_{\{h(\mathcal{T}^i) \geq \varepsilon\}}.$$

The measures $(\ell_t)_{t\geq 0}$ can, in fact, be defined simultaneously in such a way that $t \to \ell_t$ is $\Theta$-a.e. cadlag for the weak topology on finite measures on $\mathcal{T}$ ([14], Theorem 4.3). Moreover, we can $N$-a.e. recover the measure $\mu_v$, defined at (2.6), from $\mathcal{T}_v$ using the local time measures. In particular, if we define a measure $\mu^{\mathcal{T}}$ by integrating the local time measures in the following way:

$$\mu^{\mathcal{T}} := \int_0^\infty \ell_t \, dt, \tag{2.9}$$

then $\mu_v = \mu^{\mathcal{T}_v}$, $N$-a.e. This demonstrates that $\mu_v$ is indeed a natural measure for the real tree $\mathcal{T}_v$, and that it depends on the underlying excursion only through the real tree that is constructed from it. Note that alternative descriptions of $(\ell_t)_{t\geq 0}$ and $\mu^{\mathcal{T}}$ in terms of Hausdorff measures are provided in [15]. Finally, it also demonstrated in [14] that the topological support of $\mu^{\mathcal{T}}$ is $\mathcal{T}$, $\Theta$-a.e.

2.2. *Spatial trees, snakes and tours.* Consider a pair $(\mathcal{T}, \phi)$, where $\mathcal{T}$ is a compact rooted real tree and $\phi$ is a continuous mapping from $\mathcal{T}$ into $\mathbb{R}^d$; we will denote the usual Euclidean metric in $\mathbb{R}^d$ by $d_E$. We say two such pairs $(\mathcal{T}, \phi)$ and $(\mathcal{T}', \phi')$ are equivalent if and only if there exists a root preserving isometry from $\mathcal{T}$ to $\mathcal{T}'$, $\pi$, say, that also satisfies $\phi = \phi' \circ \pi$. The set of equivalence classes under this relation will be denoted $\mathbb{T}_{\text{sp}}$, and elements of this set are called spatial trees. As with real trees, we will frequently identify an equivalence class with a particular element of it. We now explain how to define a metric on this space. First, we say that a correspondence between two compact rooted real trees, $\mathcal{T}$ and $\mathcal{T}'$, is a subset $\mathcal{C} \subseteq \mathcal{T} \times \mathcal{T}'$ such that for every $\sigma \in \mathcal{T}$ there exists at least one $\sigma' \in \mathcal{T}'$ such that $(\sigma, \sigma') \in \mathcal{C}$ and, conversely, for every $\sigma' \in \mathcal{T}'$ there exists at least one $\sigma \in \mathcal{T}$ such that $(\sigma, \sigma') \in \mathcal{C}$. Moreover, we assume that $(\rho, \rho') \in \mathcal{C}$. The distortion of the correspondence $\mathcal{C}$ is defined by

$$\text{dis}(\mathcal{C}) := \sup\{|d_{\mathcal{T}}(\sigma_1, \sigma_2) - d_{\mathcal{T}'}(\sigma_1', \sigma_2')| : (\sigma_1, \sigma_1'), (\sigma_2, \sigma_2') \in \mathcal{C}\}.$$

Now define, for $(\mathcal{T}, \phi), (\mathcal{T}', \phi') \in \mathbb{T}_{\text{sp}}$, a distance by

$$d_{\text{sp}}((\mathcal{T}, \phi), (\mathcal{T}', \phi')) = \inf_{\mathcal{C} \in \mathfrak{C}(\mathcal{T}, \mathcal{T}')} \left\{ \text{dis}(\mathcal{C}) + \sup_{(\sigma, \sigma') \in \mathcal{C}} d_E(\phi(\sigma), \phi'(\sigma')) \right\}, \tag{2.10}$$

where the set $\mathfrak{C}(\mathcal{T}, \mathcal{T}')$ is the collection of all correspondences between $\mathcal{T}$ and $\mathcal{T}'$. From [14], we have that $(\mathbb{T}_{\text{sp}}, d_{\text{sp}})$ is a separable metric space (we note that it is not complete as claimed in [14]). We set $\mathcal{S} := \phi(\mathcal{T})$, which is well-defined on equivalence classes of spatial trees. Note that, although the notation $(\mathcal{T}, \phi)$ is used as shorthand for $(\mathcal{T}, d_{\mathcal{T}}, \rho, \phi)$, the notation $\mathcal{S}$ will only ever be used to denote the compact subset of $\mathbb{R}^d$ given by $\phi(\mathcal{T})$.

Moreover, when we consider the usual Hausdorff metric on compact subsets of $\mathbb{R}^d$, which we denote by $d_\mathrm{H}$, it is easy to check that if $(\mathcal{T}, \phi), (\mathcal{T}', \phi') \in \mathbb{T}_\mathrm{sp}$, then $d_\mathrm{H}(\phi(\mathcal{T}), \phi'(\mathcal{T}')) \leq d_\mathrm{sp}((\mathcal{T}, \phi), (\mathcal{T}', \phi'))$. Thus, the map from $(\mathcal{T}, \phi)$ to the compact subset $\mathcal{S}$ is continuous, and therefore measurable.

In the remainder of this section we introduce the class of spatial trees which are obtained when the mapping $\phi$ is a "Brownian embedding" of a real tree into Euclidean space, so that an arc of length $t$ in the real tree is mapped to the range of a Brownian motion run for a time $t$. Fix $x \in \mathbb{R}^d$, and let $\mathcal{T}$ be a compact rooted real tree. Consider the $\mathbb{R}^d$-valued Gaussian process $(\phi(\sigma))_{\sigma \in \mathcal{T}}$, built on a probability space with probability measure $\mathbf{P}$, whose distribution is characterized by $\mathbf{E}\phi(\sigma) = x$, $\mathrm{cov}(\phi(\sigma_1), \phi(\sigma_2)) = d_\mathcal{T}(\rho, b^\mathcal{T}(\rho, \sigma_1, \sigma_2))I$, where $I$ is the $d$-dimensional identity matrix. As remarked in [14], it is possible to chose this process to be continuous $\mathbf{P}$-a.s. for $\Theta$-a.e. $\mathcal{T}$. Assuming that we have a real tree $\mathcal{T}$ that allows us to construct a $\mathbf{P}$-a.s. continuous $\phi$, we will denote the law of $(\mathcal{T}, \phi)$ on $\mathbb{T}_\mathrm{sp}$ by $Q_x^\mathcal{T}$. This allows us to construct a $\sigma$-finite measure on $\mathbb{T}_\mathrm{sp}$ by setting

$$M_x := \int_\mathbb{T} \Theta(d\mathcal{T}) Q_x^\mathcal{T};$$

we note that the measures $Q_x^\mathcal{T}$ satisfy the necessary measurability for this integral to be well-defined. We also define $M := M_0$ and $Q^\mathcal{T} := Q_0^\mathcal{T}$. Observe that, for $M_x$-a.e. spatial tree, we have $\phi(\rho) = x$. A spatial tree chosen according to $M_x$ will be called a Brownian spatial tree started from $x$.

Pushing forward, the local time measures $\ell_t$ from $\mathcal{T}$ onto $\mathbb{R}^d$ using the map $\phi$ provides us with a cadlag (with respect to the topology induced by the weak convergence of measures on $\mathbb{R}^d$) measure-valued process $Z = (Z_t)_{t \geq 0}$. In particular, we set

$$Z_t := \ell_t \circ \phi^{-1}, \tag{2.11}$$

which defines the process $Z$ at least $M_x$-a.e. for any $x \in \mathbb{R}^d$. We will describe in Section 3 how a certain Poisson sum of these processes yields a superprocess in $\mathbb{R}^d$. The property of spatial trees that allows this connection to be made is their Markovian branching under the measures of the form $M_x$. To describe this precisely, first recall the notation $\mathrm{tr}_t(\mathcal{T})$ and $(\mathcal{T}^i)_{i \in \mathcal{I}_t}$ introduced at (2.3). The information about the spatial tree $(\mathcal{T}, \phi)$ below level $t$ is $(\mathrm{tr}_t(\mathcal{T}), \phi|_{\mathrm{tr}_t(\mathcal{T})})$, and we will denote this by $\mathcal{E}_t$. We will also write $\phi^i := \phi|_{\mathcal{T}^i}$ and $\mathcal{S}^i := \phi^i(\mathcal{T}^i)$. The Markov branching property of Brownian spatial trees can be stated as follows.

LEMMA 2.3. *Fix $t > 0$. Under the probability measure $M(\cdot | h(\mathcal{T}) > t)$ and conditional on $\mathcal{E}_t$, the collection $(\mathcal{T}^i, \phi^i)$, $i \in \mathcal{I}_t$, forms a Poisson point process on $\mathbb{T}_\mathrm{sp}$ with intensity measure*

$$\int_{\mathcal{T}_t} \ell_t(d\sigma) M_{\phi(\sigma)}. \tag{2.12}$$

PROOF. This result can be proved by a simple modification of the proof of [13], Proposition 4.2.3, using the Markov branching property of $\mathcal{T}$ that is proved in [14], Theorem 4.2. □

We define a Borel measure on $\mathcal{S}$ by setting $\mu^{\mathcal{S}} := \mu^{\mathcal{T}} \circ \phi^{-1}$, which exists and has support $\mathcal{S}$, $M_x$-a.e. It is possible to deduce that $\mu^{\mathcal{S}}$ can also be represented as $\int_0^\infty Z_t \, dt$, or, as remarked in the introduction, for $d \geq 4$, it can be expressed in terms of a Hausdorff measure on $\mathcal{S}$; see [20], Theorem 6.1. Note that this final remark immediately implies that $\mu^{\mathcal{S}}$ is a measurable function of $\mathcal{S}$, at least for $d \geq 4$ (the measurability of Hausdorff measures as functions of compact subsets of $\mathbb{R}^d$ is investigated in [25]).

In the previous section we saw how continuous excursions are useful for encoding a certain class of real trees. To perform the role of encoding the Brownian spatial trees introduced above, we will use objects called snakes and tours ([21] is a good primer for the Brownian snake). A snake is a pair of functions, $(v, w)$, say, with $v \in \mathcal{V}$ and $w$ a continuous function from $\mathbb{R}_+$ to the space of continuous paths in $\mathbb{R}^d$, which satisfy $w(s)(t) = w(s)(v(s))$ for every $t \geq v(s)$, and also $w(s)(t) = w(s')(t)$ for every $t \leq m_v(s, s')$. In fact, Theorem 2.1 of [24] shows that a snake carries some redundant information, and it suffices to consider the space of tours, where a tour is a pair $(v, r) \in C(\mathbb{R}_+, \mathbb{R}_+) \times C(\mathbb{R}_+, \mathbb{R}^d)$, with $v \in \mathcal{V}$, and $r$ a continuous $\mathbb{R}^d$-valued function which is constant on the equivalence classes given by the relation at (2.5). More precisely, in [24] it is shown that the natural map from snakes to tours given by, for $t \geq 0$, $(v(t), r(t)) = (v(t), w(t)(v(t)))$, is a homeomorphism. Due to this relationship, the process $r$ is known as the head of the snake.

The connection of snakes and tours with spatial trees can be explained as follows. Fix $(v, (\mathcal{T}_v, \phi)) \in \mathcal{V} \times \mathbb{T}_{\mathrm{sp}}$, and define, for $s \leq v(t)$, $t \geq 0$,

$$(2.13) \qquad w(t)(s) := \phi(\gamma^{\mathcal{T}_v}_{\rho,[t]}(s)),$$

where $[t]$ is the equivalence class containing $t$ under the equivalence defined at (2.5). Also define $w(t)(s) = w(t)(v(t))$ for $s \geq v(t)$. The pair $(v, w)$ is then a snake and, moreover, it is possible to check using (2.13) that the corresponding tour $(v, r)$ satisfies $r(t) = \phi([t])$. Clearly, we can recover the spatial tree from the tour $(v, r)$ by setting $\mathcal{T}_v$ to be the real tree associated with $v$ and using the final observation of the previous sentence to determine $\phi$. We can use this relationship to show that if two tours are close with respect to the uniform norm on $C(\mathbb{R}_+, \mathbb{R}_+) \times C(\mathbb{R}_+, \mathbb{R}^d)$, then so are the related spatial trees with respect to the metric $d_{\mathrm{sp}}$. The proof of this result is similar to that of Lemma 2.3 of [14], and the result obviously implies that the map $(v, r) \mapsto (\mathcal{T}, \phi)$ is measurable.

PROPOSITION 2.4. *If $(\mathcal{T}, \phi)$ and $(\mathcal{T}', \phi')$ are the spatial trees corresponding to the tours $(v, r)$ and $(v', r')$ respectively, then $d_{\mathrm{sp}}((\mathcal{T}, \phi), (\mathcal{T}', \phi')) \leq 4\|v - v'\|_\infty + \|r - r'\|_\infty$.*

PROOF. Define a correspondence between $\mathcal{T}$ and $\mathcal{T}'$ by $\mathcal{C} := \{([t], [t]') : t \in \mathbb{R}_+\}$, where $[t]$ is the equivalence class containing $t$ under the equivalence defined at (2.5) for $v$, and $[t]'$ is the corresponding quantity for $v'$. As in the proof of [14], we have $\text{dis}(\mathcal{C}) \leq 4\|v - v'\|_\infty$. It is also easy to check that $\sup_{(\sigma,\sigma') \in \mathcal{C}} d_E(\phi(\sigma), \phi'(\sigma')) = \|r - r'\|_\infty$. The proof follows on recalling the definition of $d_{\text{sp}}$ from (2.10). □

To complete this section, we will introduce the law of the Brownian tour. First, let $(v, (\mathcal{T}_v, \phi))$ be chosen in $C(\mathbb{R}_+, \mathbb{R}_+) \times \mathbb{T}_{\text{sp}}$, so that $v$ has law $N$ and, conditional on $v$, the pair $(\mathcal{T}_v, \phi)$ is a spatial tree with law $Q_x^{\mathcal{T}_v}$. Note that the marginals of this distribution are $N$ and $M_x$, so $v$ is a Brownian excursion and $(\mathcal{T}_v, \phi)$ is a Brownian spatial tree started from $x$. Define $(v, w)$ and $(v, r)$ from $(v, (\mathcal{T}_v, \phi))$ as above, and denote their laws on the appropriate spaces by $\tilde{M}'_x$ and $\tilde{M}_x$ respectively. As observed in [22], Section 6, the measure $\tilde{M}'_x$ is the law of the Brownian snake started from $x$, and we will call $\tilde{M}_x$ the law of the Brownian tour started from $x$. Accordingly, the pairs $(v, w)$ and $(v, r)$ are called the Brownian snake and Brownian tour, respectively. Note that if the subscript $x$ is missing from one of the measures defined in this paragraph, then we are working under the assumption that $x = 0$.

2.3. *CRT, ISE and normalized Brownian tour.* There is a normalization of real and spatial trees that will be of particular interest in the sections of this article where we investigate the scaling limit of simple random walks on random graph trees embedded into Euclidean space, and this is when we condition the measure $\mu^\mathcal{T}$, as defined at (2.9), to have total mass equal to one. In particular, let $N^{(1)} := N(\cdot | \tau(f) = 1)$ be the probability measure on the space of excursions $\mathcal{V}$ that is the law of the Brownian excursion, scaled to return to zero for the first time at time one. Define $\Theta^{(1)}$ from $N^{(1)}$ analogously to (2.7), and set

$$M_x^{(1)} := \int_\mathbb{T} \Theta^{(1)}(d\mathcal{T}) Q_x^\mathcal{T},$$

which is a probability measure on $\mathbb{T}_{\text{sp}}$. Also denote $M^{(1)} := M_0^{(1)}$. If $\mathcal{T}$ is a random element of $\mathbb{T}$ with law $\Theta^{(1)}$, then it is precisely (up to a deterministic scaling constant) the continuum random tree of Aldous; see [1]. Moreover, if $(\mathcal{T}, \phi)$ is a random spatial tree with law $M^{(1)}$, then we call the measure $\mu^\mathcal{S} := \mu^\mathcal{T} \circ \phi^{-1}$, which has support $\mathcal{S}$, the integrated super-Brownian excursion. This measure was first discussed by Aldous in [4], and in Section 8 we provide another characterization of it as a scaling limit, which shows how our definition matches Aldous. For an overview of the occurrence of the integrated super-Brownian excursion as the scaling limit of measures that arise in statistical mechanical models, see [27].

The normalized Brownian tour is defined similarly to the Brownian tour, except we use the above normalization in choosing the excursion. In particular, define $(v, (\mathcal{T}_v, \phi))$ to be a random variable taking values in $C([0,1], \mathbb{R}_+) \times \mathbb{T}_{\mathrm{sp}}$ such that $v$ has law $N^{(1)}$ and, conditional on $v$, the pair $(\mathcal{T}_v, \phi)$ is a spatial tree with law $Q_x^{\mathcal{T}_v}$. Defining the normalized head process $r \in C([0,1], \mathbb{R}^d)$ similarly to above, we call $(v, r)$ the normalized Brownian tour started from $x$, and its law will be denoted $\tilde{M}_x^{(1)}$, with $\tilde{M}^{(1)} := \tilde{M}_0^{(1)}$.

**3. Homeomorphism between $\mathcal{T}$ and $\mathcal{S}$.** The purpose of this section is to show that when $d \geq 8$ the continuous map $\phi$ is actually a homeomorphism from $(\mathcal{T}, d_{\mathcal{T}})$ to $(\mathcal{S}, d_E)$, $M$-a.e. In fact, because $(\mathcal{T}, d_{\mathcal{T}})$ is compact and $(\mathcal{S}, d_E)$ is Hausdorff, $M$-a.e., it is sufficient to show that $\phi$ is injective, $M$-a.e. We will prove that this is the case by applying results proved in [11] about the multiple points of super-processes. Although our main conclusion is suggested by comments made in Section 3.4 of [4], we can find no rigorous proof in the literature and so continue to prove it here.

We start by describing briefly the connection between spatial trees and super-processes that has been developed by Duquesne and Le Gall for further details see [13] and [14]. First, let $(\mathcal{T}^i, \phi^i)_{i \in \mathcal{I}}$ be a Poisson process on $\mathbb{T}_{\mathrm{sp}}$ with intensity measure $M$. For each $i \in \mathcal{I}$, define $(Z_t^i)_{t \geq 0}$ to be the measure-valued process associated with $(\mathcal{T}^i, \phi^i)$ by the formula at (2.11). If we then define the process $Y = (Y_t)_{t \geq 0}$ by setting, for $t > 0$,

$$Y_t := \sum_{i \in \mathcal{I}} Z_t^i, \tag{3.1}$$

and $Y_0 := \delta_0$, where $\delta_0$ is the probability measure on $\mathbb{R}^d$ that places all of its mass at the origin, then $Y$ is a Dawson–Watanabe super-process, started from $\delta_0$ ([14], Proposition 6.1). We assume that the Poisson process $(\mathcal{T}^i, \phi^i)_{i \in \mathcal{I}}$ is built on an underlying probability space with probability measure $\mathbf{P}$.

In particular, it will be important for the arguments that we apply to have a description of the range of the super-process $Y$ in terms of the sets $(\mathcal{S}_t^i)_{t > 0, i \in \mathcal{I}}$, where $\mathcal{S}_t^i := \phi^i(\mathcal{T}_t^i)$, and $\mathcal{T}_t^i$ is the level $t$ subset of $\mathcal{T}^i$ defined by (2.2). For a Borel measure $\nu$ on $\mathbb{R}^d$, denote by $S(\nu)$ it closed support in $\mathbb{R}^d$, and set, for $s \leq t$, $\mathcal{R}(s,t)$ to be the closure of $\bigcup_{s \leq r \leq t} S(Y_r)$ with respect to the Euclidean metric. Furthermore, define

$$\begin{aligned}\mathcal{R}_1 &:= \bigcup_{0 < s \leq t < \infty} \mathcal{R}(s,t), \\ \mathcal{R}_2 &:= \bigcup_{0 < s_1 \leq t_1 < s_2 \leq t_2 < \infty} \mathcal{R}(s_1, t_1) \cap \mathcal{R}(s_2, t_2),\end{aligned} \tag{3.2}$$

which are the range (or one-multiple points) and two-multiple points of $Y$ respectively. As noted in the proof of [14], Proposition 6.2, $\mathbf{P}$-a.s., we have that $S(Y_t) \subseteq \bigcup_{i \in \mathcal{I}} \mathcal{S}_t^i$ for every $t > 0$, with equality holding for all but a countable collection of times, $D$, say. Furthermore, since $M(h(\mathcal{T}) > s)$ is finite for any $s > 0$ and $\phi^i$ is continuous for every $i \in \mathcal{I}$, **P**-a.s., there can only be a finite number of sets $\bigcup_{s \leq r \leq t} \mathcal{S}_r^i$ which are nonempty, and because each set of the form $\bigcup_{s \leq r \leq t} \mathcal{S}_r^i$ is closed, we must have that $\mathcal{R}(s,t) \subseteq \bigcup_{s \leq r \leq t} \bigcup_{i \in \mathcal{I}} \mathcal{S}_r^i$ for every $s \leq t$, **P**-a.s. Using the continuity of the maps $\phi^i$ and the countability of $D$, we are also able to deduce that $\bigcup_{s \leq r \leq t} \bigcup_{i \in \mathcal{I}} \mathcal{S}_r^i \subseteq \mathcal{R}(s,t)$ for every $s \leq t$, $s \notin D$, **P**-a.s. Hence, again applying the countability of $D$, we have the following alternative expression for $\mathcal{R}_1$, **P**-a.s.:

$$\mathcal{R}_1 = \bigcup_{t > 0} \bigcup_{i \in \mathcal{I}} \mathcal{S}_t^i. \tag{3.3}$$

Before we proceed with our main argument, we collect some other properties of $\mathcal{R}_1$ and $\mathcal{R}_2$ that are proved in [11], Theorems 1.5 and 1.6. With regards to part (a) of the following lemma, note that the precise Hausdorff measure functions of $\mathcal{R}_1$ for $d \geq 4$ are found in [20], Theorem 1.1.

LEMMA 3.1. (a) *If $d > 4$, $\mathcal{R}_1$ has $\sigma$-finite Hausdorff measure with respect to the measure function $x^4 \ln \ln x^{-1}$, **P**-a.s.*

(b) *If $d > 4$ and $A \subseteq \mathbb{R}^d$ is null for the Hausdorff measure with respect to the measure function $x^{d-4}$, then $\mathcal{R}_1 \cap A = \varnothing$, **P**-a.s.*

(c) *If $d \geq 8$, then $\mathcal{R}_2 = \varnothing$, **P**-a.s.*

A useful corollary of parts (a) and (b) of this result is the following, which demonstrates that "independent" spatial trees do not intersect when they are started from different points of $\mathbb{R}^d$ and $d$ is large.

COROLLARY 3.2. (a) *Fix $x, y \in \mathbb{R}^d$. If $d \geq 8$, then $M_x \otimes M_y$-a.e. we have $\mathcal{S} \cap \mathcal{S}' = \{x\} \cap \{y\}$, where $\mathcal{S} := \phi(\mathcal{T})$, $\mathcal{S}' := \phi'(\mathcal{T}')$ and $((\mathcal{T}, \phi), (\mathcal{T}', \phi'))$ represents an element of $\mathbb{T}_{\rm sp}^2$.*

(b) *Part (a) holds when the measure $M_x \otimes M_y$ is replaced by the probability measure $M_x(\cdot | h(\mathcal{T}) > t) \otimes M_y(\cdot | h(\mathcal{T}) > t)$, for any $t > 0$.*

PROOF. By Lemma 3.1(a), we know that $\mathcal{R}_1$ has $\sigma$-finite Hausdorff measure with respect to the measure function $x^4 \ln \ln x^{-1}$, **P**-a.s. Thus, from the expression at (3.3) and the Poisson process construction of the super-process, we see that $\mathcal{S} = \phi(\mathcal{T})$ also satisfies this property, $M$-a.s. Since we are assuming that $d \geq 8$, it follows that $\mathcal{S}$ is null for $(d-4)$-dimensional Hausdorff measure, $M_x$-a.s.

Now suppose that $(\mathcal{S}, y + \mathcal{R}_1)$ is chosen according to $M_x \otimes \mathbf{P}$, where $\mathcal{R}_1$ is the range of the super-process as described above. From Lemma 3.1(b) and the conclusion of the previous paragraph, we have that $(y + \mathcal{R}_1) \cap \mathcal{S} = \varnothing$, $M_x \otimes \mathbf{P}$-a.s. Again applying (3.3) and the Poisson process description of the super-process, it follows that $\mathcal{S} \cap \mathcal{S}' \subseteq \{y\}$, $M_x \otimes M_y$-a.s. Using symmetry, part (a) is a straightforward consequence of this. Part (b) follows immediately. $\square$

In the next lemma we combine the above result with the Markov branching property of spatial trees to deduce the disjointness of a particular collection of subsets of spatial trees. Recall from Section 2.1 the definition of subtrees above a certain level, $(\mathcal{T}^i)_{i \in \mathcal{I}_t}$, and also the definition of the spatial trees $(\mathcal{T}^i, \phi^i)_{i \in \mathcal{I}_t}$ from Section 2.2.

LEMMA 3.3. *Let $d \geq 8$ and fix $t > 0$. Under the measure $M$, the sets $\mathcal{S}^{i,o}$, $i \in \mathcal{I}_t$, where $\mathcal{S}^{i,o} := \phi^i(\mathcal{T}^i) \setminus \{\phi^i(\rho^i)\}$, are almost surely disjoint.*

PROOF. Write $M' := M(\cdot | h(\mathcal{T}) > t)$. From Lemma 2.3 we have that under the probability measure $M'$, conditional on $\mathcal{E}_t$, the collection $(\mathcal{T}^i, \phi^i)$, $i \in \mathcal{I}_t$ is a Poisson process on $\mathbb{T}_{\text{sp}}$ with intensity measure given by (2.12). Now, choose $\varepsilon > 0$ and note that it is possible to deduce from (2.8) and Theorem 2.2(b) that $\int_{\mathcal{T}_t} \ell_t(d\sigma) M_{\phi(\sigma)}(h(\mathcal{T}) > \varepsilon) \in (0, \infty)$, $M'$-a.e. Thus, it makes sense to further condition on the size of $\mathcal{I}_{t,\varepsilon} := \{i \in \mathcal{I}_t : h(\mathcal{T}^i) > \varepsilon\}$. In particular, under the measure $M'(\cdot | \mathcal{E}_t, \#\mathcal{I}_{t,\varepsilon} = n)$, the elements of $\{(\mathcal{T}^i, \phi^i) : i \in \mathcal{I}_{t,\varepsilon}\}$ are distributed as a sample of $n$ independent random variables, each with law

$$\int_{\mathcal{T}_t} \ell_t(d\sigma) M_{\phi(\sigma)}(\cdot | h(\mathcal{T}) > \varepsilon) \ell_t(\mathcal{T}_t).$$

By first conditioning on the locations of the points $\phi^i(\rho^i)$, it is straightforward to apply Corollary 3.2(b) to deduce from this that

$$M'((\mathcal{S}^{i,o})_{i \in \mathcal{I}_{t,\varepsilon}} \text{ are disjoint}) = 1.$$

There is no problem in extending this result to deduce that

$$M'((\mathcal{S}^{i,o})_{i \in \mathcal{I}_t} \text{ are disjoint}) = 1,$$

which completes the proof, because on $\{h(\mathcal{T}) \leq t\}$ the set $\mathcal{I}_t$ is empty. $\square$

The above result provides the first ingredient in our proof of the fact that $\phi$ is injective. The second is given by the following lemma, which shows that the image under $\phi$ of the level sets of $\mathcal{T}$ are disjoint.

LEMMA 3.4. *Let $d \geq 8$. $M$-a.e., the sets $\mathcal{S}_t := \phi(\mathcal{T}_t)$, $t \geq 0$, are disjoint.*

PROOF. First assume that at all but a countable collection of points, $D$, say, we have $S(Y_t) = \bigcup_{i \in \mathcal{I}} \mathcal{S}_t^i$, where $Y$ is the measure-valued process defined at (3.1) from the Poisson collection of spatial trees. As remarked earlier in this section, a proof that this fact holds $\mathbf{P}$-a.s. appears within the proof of [14], Proposition 6.2.

Now suppose there exists an $x \in \mathbb{R}^d$ such that $x \in (\bigcup_{i \in \mathcal{I}} \mathcal{S}_s^i) \cap (\bigcup_{i \in \mathcal{I}} \mathcal{S}_t^i)$, for some $0 < s < t$. Set $\varepsilon := (t-s)/2$. Clearly, for some $i \in \mathcal{I}$, we can find $\sigma \in \mathcal{T}^i$ such that $\phi^i(\sigma) = x$ and $d_{\mathcal{T}^i}(\rho^i, \sigma) = s$, where $\rho^i$ is the root of $\mathcal{T}^i$. By considering the arc in $\mathcal{T}^i$ from $\rho^i$ to $\sigma$, we can find a sequence $(\sigma_n)_{n \geq 0}$ that converges in $\mathcal{T}^i$ to $\sigma$ with $d_{\mathcal{T}^i}(\rho^i, \sigma_n) \in [0, s] \setminus D$ for each $n$. Similarly, for some $j \in \mathcal{I}$, we can find $\sigma' \in \mathcal{T}^j$ such that $\phi^j(\sigma') = x$ and $d_{\mathcal{T}^j}(\rho^j, \sigma') = t$, and also a sequence $(\sigma'_n)_{n \geq 0}$ converging in $\mathcal{T}^j$ to $\sigma'$ with $d_{\mathcal{T}^j}(\rho^j, \sigma'_n) \in [0, t] \setminus D$ for each $n$. It follows that there exists an $N$ such that, for $n \geq N$,

$$d_{\mathcal{T}^i}(\rho^i, \sigma_n) \in [s/2, s], \qquad d_{\mathcal{T}^j}(\rho^j, \sigma'_n) \in [t-\varepsilon, t].$$

Since $d_{\mathcal{T}^i}(\rho^i, \sigma_n), d_{\mathcal{T}^j}(\rho^j, \sigma'_n) \notin D$, we have, for $n \geq N$,

$$\phi^i(\sigma_n) \in \mathcal{S}^i_{d_{\mathcal{T}^i}(\rho^i, \sigma_n)} \subseteq S(Y_{d_{\mathcal{T}^i}(\rho^i, \sigma_n)}) \subseteq \mathcal{R}(s/2, s),$$

and also $\phi^j(\sigma'_n) \in \mathcal{R}(t-\varepsilon, t)$. By continuity, we have $\phi^i(\sigma_n), \phi^j(\sigma'_n) \to x$. Hence, $x \in \mathcal{R}(s/2, s) \cap \mathcal{R}(t-\varepsilon, t) \subseteq \mathcal{R}_2$, where $\mathcal{R}_2$ is the set of two-multiple points of our super-process. However, by Lemma 3.1(c), $\mathcal{R}_2 = \varnothing$, and so no such $x$ exists, $\mathbf{P}$-a.s. Consequently, the sets $\bigcup_{i \in \mathcal{I}} \mathcal{S}_t^i, t > 0$ are disjoint, $\mathbf{P}$-a.s.

Observe that $\bigcup_{i \in \mathcal{I}} \mathcal{S}_0^i = \{0\}$ and also, from (3.3), $\bigcup_{t > 0} \bigcup_{i \in \mathcal{I}} \mathcal{S}_t^i = \mathcal{R}_1$. Hence, to show that the sets $\bigcup_{i \in \mathcal{I}} \mathcal{S}_t^i, t \geq 0$ are disjoint, $\mathbf{P}$-a.s., it suffices to show that $\mathcal{R}_1 \cap \{0\} = \varnothing$, $\mathbf{P}$-a.s. However, this is a consequence of Lemma 3.1(b). The result now follows easily on recalling the Poisson process definition of the super-process $Y$. □

We are now ready to proceed with the main result of this section.

PROPOSITION 3.5. *Let $d \geq 8$. $M$-a.e., the map $\phi : \mathcal{T} \to \mathcal{S}$ is injective.*

PROOF. The following proof holds $M$-a.e. Suppose that there exist distinct $\sigma_1, \sigma_2 \in \mathcal{T}$ that satisfy $\phi(\sigma_1) = \phi(\sigma_2)$. By the previous lemma, we know that $\sigma_1, \sigma_2 \in \mathcal{T}_t$, for some $t > 0$. Necessarily we must also have that $b^{\mathcal{T}}(\rho, \sigma_1, \sigma_2) \in \mathcal{T}_s$ for some $s < t$ [recall the notation for a branch-point of $\mathcal{T}$ from (2.1)]. Choose $r \in (s, t) \cap \mathbb{Q}$, and consider the collection of spatial trees above level $r$, which, using the notation introduced above Lemma 2.3, can be written as $(\mathcal{T}^i, \phi^i), i \in \mathcal{I}_r$. Since $b^{\mathcal{T}}(\rho, \sigma_1, \sigma_2) \in \mathcal{T}_s$ for some $s < r$, we have that $\sigma_1 \in \mathcal{T}^i$ and $\sigma_2 \in \mathcal{T}^j$ for some $i \neq j$. Moreover, if $\phi(\sigma_1) = \phi^i(\rho^i)$, then $\phi(\sigma_1) \in \mathcal{S}_r \cap \mathcal{S}_t$. However, by Lemma 3.4, we can assume that $\mathcal{S}_r \cap \mathcal{S}_t = \varnothing$, and so $\phi(\sigma_1) \in \mathcal{S}^{i,o}$. Thus, by symmetry, we have that $\phi(\sigma_1) \in \mathcal{S}^{i,o} \cap \mathcal{S}^{j,o}$.

After extending Lemma 3.3 using a countability argument, the sets $\mathcal{S}^{i,o}$, $i \in \mathcal{I}_q$, can be assumed to be disjoint for any $q \in \mathbb{Q}$. In particular, the sets $\mathcal{S}^{i,o}$, $i \in \mathcal{I}_r$, are disjoint. Consequently, no such $\sigma_1$ and $\sigma_2$ exist, which implies that $\phi$ is injective. $\square$

**4. Hausdorff measure of spatial tree arcs.** In the previous section we showed that when the dimension $d$ is large enough, the map $\phi$ is a homeomorphism, $M$-a.e. Consequently, the set $\mathcal{S}$ is a dendrite and, as remarked in the Introduction, between any two points of $\mathcal{S}$ there is a unique arc in $\mathcal{S}$, $M$-a.e. We show in this section how the natural way to measure the distance along the arcs of $\mathcal{S}$ is to use the Hausdorff measure with respect to the measure function $c_d x^2 \ln \ln x^{-1}$, where $c_d$ is a deterministic constant that depends only on the dimension $d$. The following result will be fundamental in proving this; it determines the Hausdorff measure of Brownian paths in high dimensions. The description of the Hausdorff measure $\mathcal{H}$ that appears in the lemma should be considered to be a definition.

LEMMA 4.1 ([7], Theorem 5). *Suppose that $(B_t)_{t \geq 0}$ is a standard $d$-dimensional Brownian motion built on a probability space $(\Omega, \mathcal{F}, \mathbf{P})$. If $d \geq 3$, then*

$$\mathcal{H}(\{B_s : s \in [0,t]\}) = t \qquad \forall t \geq 0, \mathbf{P}\text{-}a.s.,$$

*where $\mathcal{H}$ is the Hausdorff measure calculated with respect to the function $c_d x^2 \ln \ln x^{-1}$, and $c_d$ is a deterministic constant that depends only upon $d$.*

Recall from Section 2.1 that the path in $\mathcal{T}$ of unit speed between $\sigma_1$ and $\sigma_2$ is denoted by $\gamma^{\mathcal{T}}_{\sigma_1,\sigma_2}$, and its image by $\Gamma^{\mathcal{T}}_{\sigma_1,\sigma_2}$. The following result applies the above lemma to describe the Hausdorff measure of the sets $\phi(\Gamma^{\mathcal{T}}_{\rho,\sigma})$ for $\sigma \in \mathcal{T}$.

LEMMA 4.2. *Let $d \geq 3$. For $M$-a.e. $(\mathcal{T}, \phi)$,*

$$\mathcal{H}(\phi(\Gamma^{\mathcal{T}}_{\rho,\sigma})) = d_{\mathcal{T}}(\rho, \sigma) \qquad \forall \sigma \in \mathcal{T}.$$

PROOF. Fix a countably dense sequence $(t_n^*)_{n \geq 0}$ in $\mathbb{R}_+$. Given an excursion $v \in \mathcal{V}$, we will denote by $\sigma_n^*$ the vertex $[t_n^*]$ of the corresponding real tree $\mathcal{T}$, where $[t]$ represents the equivalence class of $t$ under the equivalence defined at (2.5).

We start the proof of the lemma by demonstrating the claim that $\tilde{M}$-a.e.,

(4.1) $$\mathcal{H}(\phi(\Gamma^{\mathcal{T}}_{\rho,\sigma})) = d_{\mathcal{T}}(\rho, \sigma) \qquad \forall \sigma \in \Gamma^{\mathcal{T}}_{\rho,\sigma_n^*}, n \geq 0,$$

where $\tilde{M}$ is the law of the Brownian tour introduced in Section 2.2. By countability, it will suffice to prove the above result holds for one particular $n$. Observe that, conditional on $v$, the process $\phi \circ \gamma^{\mathcal{T}}_{\rho,\sigma_n^*}$ is a standard $d$-dimensional Brownian motion run for a time $v(t_n^*) = d_\mathcal{T}(\rho, \sigma_n^*)$. Consequently, by Lemma 4.1, it satisfies

$$(4.2) \qquad \mathcal{H}(\{\phi(\gamma_{\rho,\sigma_n^*}^\mathcal{T}(s)) : 0 \leq s \leq t\}) = t \qquad \forall t \in [0, d_\mathcal{T}(\rho, \sigma_n^*)],$$

$\tilde{M}$-a.e. Now, for $\sigma \in \Gamma_{\rho,\sigma_n^*}^\mathcal{T}$, we have that $\sigma = \gamma_{\rho,\sigma_n^*}^\mathcal{T}(d_\mathcal{T}(\rho,\sigma))$, because $\gamma_{\rho,\sigma_n^*}^\mathcal{T}$ traverses the arc $\Gamma_{\rho,\sigma_n^*}^\mathcal{T}$ at a unit speed. Furthermore, it is clear from the definitions that $\gamma_{\rho,\sigma_n^*}^\mathcal{T}|_{[0,d_\mathcal{T}(\rho,\sigma)]} = \gamma_{\rho,\sigma}^\mathcal{T}$ and, thus, $\gamma_{\rho,\sigma_n^*}^\mathcal{T}([0,d_\mathcal{T}(\rho,\sigma)]) = \Gamma_{\rho,\sigma}^\mathcal{T}$. Hence, we can rewrite (4.2) to obtain that $\mathcal{H}(\phi(\Gamma_{\rho,\sigma}^\mathcal{T})) = d_\mathcal{T}(\rho,\sigma)$ for $\sigma \in \Gamma_{\rho,\sigma_n^*}^\mathcal{T}$, $\tilde{M}$-a.e., which completes the proof of the claim at (4.1).

Suppose we have a realization of $(v, r)$ for which the claim at (4.1) holds, and let $\sigma \in \mathcal{T}$. Since $\mathcal{T}^* := \{\sigma_n^* : n \geq 0\}$ is dense in $\mathcal{T}$, there exists a sequence $(\sigma_n)_{n \geq 0}$ in $\mathcal{T}^*$ such that $d_\mathcal{T}(\sigma_n, \sigma) \to 0$. Clearly, if we define $b_n := b^\mathcal{T}(\rho, \sigma_n, \sigma)$, then $b_n \in \Gamma_{\rho,\sigma_n}^\mathcal{T}$ for each $n$, and so the claim at (4.1) implies that $\mathcal{H}(\phi(\Gamma_{\rho,b_n}^\mathcal{T})) = d_\mathcal{T}(\rho, b_n)$. It is also straightforward to check that $d_\mathcal{T}(b_n, \sigma) \to 0$, and so

$$(4.3) \qquad \lim_{n \to 0} \mathcal{H}(\phi(\Gamma_{\rho,b_n}^\mathcal{T})) = d_\mathcal{T}(\rho, \sigma).$$

Furthermore, since $d_\mathcal{T}(\rho, b_n) \leq d_\mathcal{T}(\rho, \sigma)$ for every $n$, it is possible to choose a subsequence of $(b_{n_i})_{i \geq 0}$ such that $d_\mathcal{T}(\rho, b_{n_i})$ is increasing. It follows that the set sequence $(\Gamma_{\rho,b_{n_i}}^\mathcal{T})_{i \geq 0}$ is also increasing and we must have that $\bigcup_i \Gamma_{\rho,b_{n_i}}^\mathcal{T}$ is equal to either $\Gamma_{\rho,\sigma}^\mathcal{T}$ or $\Gamma_{\rho,\sigma}^\mathcal{T} \setminus \{\sigma\}$. Since the $\mathcal{H}$-measure of a set is unaffected by removing one point, we have that $\mathcal{H}(\phi(\Gamma_{\rho,\sigma}^\mathcal{T})) = \lim_{i \to \infty} \mathcal{H}(\phi(\Gamma_{\rho,b_{n_i}}^\mathcal{T}))$, and combining this with the result at (4.3) implies that $\mathcal{H}(\phi(\Gamma_{\rho,\sigma}^\mathcal{T})) = d_\mathcal{T}(\rho, \sigma)$.

From the conclusions of the two previous paragraphs we obtain that $\mathcal{H}(\phi(\Gamma_{\rho,\sigma}^\mathcal{T})) = d_\mathcal{T}(\rho, \sigma)$ for every $\sigma \in \mathcal{T}$, $\tilde{M}$-a.e., and since the marginal law of $(\mathcal{T}, \phi)$ under $\tilde{M}$ is $M$, this completes the proof. $\square$

As remarked above, when $\phi$ is injective, $\mathcal{S}$ is a dendrite, and there exists a unique arc between any two points in $\mathcal{S}$. We will denote the arc between $x_1$ and $x_2$ in $\mathcal{S}$ by $\Gamma_{x_1,x_2}^\mathcal{S}$; this is defined to be the image of any continuous injection from $[0,1]$ to $\mathcal{S}$ that takes the value $x_1$ at zero and $x_2$ at one. Observe that when $\phi$ is injective, we clearly have $\Gamma_{\phi(\sigma_1),\phi(\sigma_2)}^\mathcal{S} = \phi(\Gamma_{\sigma_1,\sigma_2}^\mathcal{T})$, for any $\sigma_1, \sigma_2 \in \mathcal{T}$. The main result of this section is the following, which describes precisely the Hausdorff measure of arcs in $\mathcal{S}$.

PROPOSITION 4.3. *Let $d \geq 8$. For $M$-a.e. $(\mathcal{T}, \phi)$,*

$$\mathcal{H}(\Gamma_{\phi(\sigma_1),\phi(\sigma_2)}^\mathcal{S}) = d_\mathcal{T}(\sigma_1, \sigma_2) \qquad \forall \sigma_1, \sigma_2 \in \mathcal{T}.$$

PROOF. By Proposition 3.5 and Lemma 4.2, we can assume that $\phi$ is injective and $\mathcal{H}(\Gamma^{\mathcal{S}}_{0,\phi(\sigma)}) = d_{\mathcal{T}}(\rho,\sigma)$ for every $\sigma \in \mathcal{T}$. Applying this and the identity
$$\Gamma^{\mathcal{S}}_{\phi(\sigma_1),\phi(\sigma_2)} = (\Gamma^{\mathcal{S}}_{0,\phi(\sigma_1)} \setminus \Gamma^{\mathcal{S}}_{0,\phi(b)}) \cup (\Gamma^{\mathcal{S}}_{0,\phi(\sigma_2)} \setminus \Gamma^{\mathcal{S}}_{0,\phi(b)}) \cup \{\phi(b)\},$$
which holds for each $\sigma_1, \sigma_2 \in \mathcal{T}$, where $b$ is the branch-point of $\rho$, $\sigma_1$ and $\sigma_2$ in $\mathcal{T}$, it is readily checked that $\mathcal{H}(\Gamma^{\mathcal{S}}_{\phi(\sigma_1),\phi(\sigma_2)}) = d_{\mathcal{T}}(\rho,\sigma_1) + d_{\mathcal{T}}(\rho,\sigma_2) - 2d_{\mathcal{T}}(\rho,b)$, from which the result follows. □

**5. Recovering spatial trees in high dimensions.** The result about the Hausdorff measure of Brownian paths stated as Lemma 4.1 can be used to recover the path $(B_t)_{t \geq 0}$ from its range, $\mathcal{R} := \{B_t : t \geq 0\}$, when $d \geq 4$. In particular, let $d \geq 4$, so that $t \mapsto B_t$ is injective and $\mathcal{R}$ is a dendrite, **P**-a.s. If we define the function $H : \mathcal{R} \to \mathbb{R}_+$ by $H(x) := \mathcal{H}(\Gamma^{\mathcal{R}}_{0,x})$, where $\Gamma^{\mathcal{R}}_{0,x}$ is the unique arc between 0 and $x$ in $\mathcal{R}$, then it is easy to check using Lemma 4.1 that $H^{-1}(t) = B_t$ for all $t \geq 0$, **P**-a.s. Thus, the following result holds.

LEMMA 5.1. *Suppose that $B = (B_t)_{t \geq 0}$ is a standard $d$-dimensional Brownian motion built on a probability space $(\Omega, \mathcal{F}, \mathbf{P})$, and $\mathcal{R} := \{B_t : t \geq 0\}$ is its range. If $d \geq 4$, then there exists a set $\Omega^* \subseteq \Omega$ such that $\mathbf{P}(\Omega \setminus \Omega^*) = 0$ and also, if $\omega, \tilde{\omega} \in \Omega^*$, then*
$$\mathcal{R}^{\omega} = \mathcal{R}^{\tilde{\omega}} \quad \Leftrightarrow \quad B^{\omega} = B^{\tilde{\omega}},$$
*where the superscript $\omega$ illustrates the dependence of the random variables on $\omega \in \Omega$.*

In this section we will show an analogous result demonstrating that it is possible to recover the spatial tree $(\mathcal{T}, \phi)$ from the compact set $\phi(\mathcal{T}) \subseteq \mathbb{R}^d$, $M$-a.e., if the dimension $d$ is large enough. As a consequence of this, we will also exhibit how to recover the super-process $Y$, as defined at (3.1), from its range $\mathcal{R}_1$. As in the case of recovering a Brownian path from its range, the key to our proof will be using the Hausdorff measure $\mathcal{H}$ to measure distance along arcs.

When $\phi$ is injective, $\mathcal{S} = \phi(\mathcal{T})$ is a dendrite, and so the function $d_{\mathcal{S}} : \mathcal{S} \times \mathcal{S} \to \mathbb{R}_+$ obtained by setting
$$d_{\mathcal{S}}(x_1, x_2) := \mathcal{H}(\Gamma^{\mathcal{S}}_{x_1,x_2})$$
is well-defined. We now show that, in fact, $(\mathcal{S}, d_{\mathcal{S}}, 0)$ is a real tree equivalent to $(\mathcal{T}, d_{\mathcal{T}}, \rho)$, $M$-a.e.

PROPOSITION 5.2. *Suppose $d \geq 8$. For $M$-a.e. choice of spatial tree $(\mathcal{T}, \phi)$, the pointed metric space $(\mathcal{S}, d_{\mathcal{S}}, 0)$ is a rooted real tree equivalent to $(\mathcal{T}, d_{\mathcal{T}}, \rho)$ and, moreover, if we define the map $I : \mathcal{S} \to \mathbb{R}^d$ to be the restriction of the identity map in $\mathbb{R}^d$, then $(\mathcal{S}, I)$ and $(\mathcal{T}, \phi)$ are equivalent spatial trees.*

PROOF. By Proposition 4.3, we can assume that $\phi$ is a bijection and $d_{\mathcal{S}}(\phi(\sigma_1), \phi(\sigma_2)) = d_{\mathcal{T}}(\sigma_1, \sigma_2)$, for every $\sigma_1, \sigma_2 \in \mathcal{T}$. Thus, $\phi : (\mathcal{T}, d_{\mathcal{T}}) \to (\mathcal{S}, d_{\mathcal{S}})$ is actually an isometry. Consequently, because we also have that $\phi(\rho) = 0$, the pointed metric space $(\mathcal{S}, d_{\mathcal{S}}, 0)$ is a rooted real tree equivalent to $(\mathcal{T}, d_{\mathcal{T}}, \rho)$. The equivalence of spatial trees is a result of the identity $\phi \equiv I \circ \phi$.
□

A corollary of this result is that, if we are given the set $\mathcal{S}$, we can determine the spatial tree $(\mathcal{T}, \phi)$ that was used to construct it, $M$-a.e. Recalling that under $M$ the tree $\mathcal{T}$ is constructed from a Brownian excursion and $\phi$ is a Brownian embedding, the following corollary is reminiscent of the result proved in [6] about recovering from an iterated Brownian motion the two underlying Brownian motions used in its construction. Also closely related to this result is Corollary 5.5.

COROLLARY 5.3. *For $d \geq 8$, there exists a set $\mathbb{T}_{\mathrm{sp}}^* \subseteq \mathbb{T}_{\mathrm{sp}}$ that satisfies $M(\mathbb{T}_{\mathrm{sp}} \setminus \mathbb{T}_{\mathrm{sp}}^*) = 0$ and also, if $(\mathcal{T}, \phi), (\tilde{\mathcal{T}}, \tilde{\phi}) \in \mathbb{T}_{\mathrm{sp}}^*$, then the compact sets $\phi(\mathcal{T})$ and $\tilde{\phi}(\tilde{\mathcal{T}})$ are equal if and only if $(\mathcal{T}, \phi)$ and $(\tilde{\mathcal{T}}, \tilde{\phi})$ are equivalent spatial trees.*

We continue by presenting a further corollary of Proposition 5.2 which shows that given the range of a certain super-process we can determine the super-process itself if the dimension is large enough.

COROLLARY 5.4. *Suppose that $Y$ is the Dawson–Watanabe super-process in $\mathbb{R}^d$, $d \geq 8$, started from $\delta_0$, with range $\mathcal{R}_1$, built on a probability space $(\Omega, \mathcal{F}, \mathbf{P})$. There exists a set $\Omega^* \subseteq \Omega$ such that $\mathbf{P}(\Omega \setminus \Omega^*) = 0$ and also, if $\omega, \tilde{\omega} \in \Omega^*$, then*

$$\mathcal{R}_1^\omega = \mathcal{R}_1^{\tilde{\omega}} \quad \Leftrightarrow \quad Y^\omega = Y^{\tilde{\omega}},$$

*where the superscript $\omega$ illustrates the dependence of the random variables on $\omega \in \Omega$.*

PROOF. As in Section 3, we can assume that $Y$ is built from a Poisson process of spatial trees $(\mathcal{T}^i, \phi^i)_{i \in \mathcal{I}}$, with intensity measure $M$. By Corollary 5.3, we can, in fact, regard this as a Poisson process on $\mathbb{T}_{\mathrm{sp}}^*$, so that each $\mathcal{S}^i := \phi^i(\mathcal{T}^i)$ determines the spatial tree $(\mathcal{T}^i, \phi^i) \in \mathbb{T}_{\mathrm{sp}}^*$ uniquely, $\mathbf{P}$-a.s. Furthermore, by applying an argument almost identical to that of Lemma 3.3, we are able to deduce that the sets $\mathcal{S}^{i,o} := \phi^i(\mathcal{T}^i) \setminus \{0\}$, $i \in \mathcal{I}$, are disjoint $\mathbf{P}$-a.s, and also, from (3.3), we have that $\mathcal{R}_1 = \bigcup_{i \in \mathcal{I}} \bigcup_{t>0} \mathcal{S}_t^i$, $\mathbf{P}$-a.s. By the injectivity of the maps $\phi^i$, which follows from Proposition 3.5, we must have that $\bigcup_{t>0} \mathcal{S}_t^i = \mathcal{S}^{i,o}$ for each $i$, $\mathbf{P}$-a.s. Combining these facts, we can conclude that there exists a set $\Omega^* \subseteq \Omega$ which has probability one and upon which $\mathcal{R}_1 = \bigcup_{i \in \mathcal{I}} \mathcal{S}^{i,o}$, $(\mathcal{S}^{i,o})_{i \in \mathcal{I}}$ are disjoint, and $(\mathcal{T}^i, \phi^i) \in \mathbb{T}^*_{\text{sp}}$ for each $i \in \mathcal{I}$.

Now consider $\mathcal{R}_1^\omega$, for some $\omega \in \Omega^*$. Applying the first two properties that are assumed to hold on $\Omega^*$, we are able to deduce that the set of path-connected components of $\mathcal{R}_1^\omega$, $\mathcal{C}^\omega$, say, is precisely equal to the set $\{\mathcal{S}^{i,o} : i \in \mathcal{I}^\omega\}$. Hence, the set $\mathcal{C}_0^\omega := \{A \cup \{0\} : A \in \mathcal{C}^\omega\}$ must be equal to $\{\mathcal{S}^i : i \in \mathcal{I}^\omega\}$. Therefore, by the definition of $\mathbb{T}^*_{\text{sp}}$, each set $A \in \mathcal{C}_0^\omega$ determines uniquely a spatial tree $(\mathcal{T}_A, \phi_A) \in \mathbb{T}^*_{\text{sp}}$. Now the collection $\{(\mathcal{T}_A, \phi_A) : A \in \mathcal{C}_0^\omega\}$ is completely determined by $\mathcal{R}_1^\omega$, and is equal to $\{(\mathcal{T}^i, \phi^i) : i \in \mathcal{I}^\omega\}$. Thus, $\mathcal{R}_1^\omega$ determines the super-process $Y^\omega$, and if $\mathcal{R}_1^\omega = \mathcal{R}_1^{\tilde\omega}$ for some $\omega, \tilde\omega \in \Omega^*$ then, $Y^\omega = Y^{\tilde\omega}$ as claimed. □

Finally, we state the analogue of Corollary 5.3 in the case of ordered trees. More precisely, the previous results of this section lead easily to the fact that in high dimensions the Brownian tour, $(v, r)$ and, consequently, the Brownian snake, are determined by the Brownian head process, $r$, thus, showing that there is enough information in the ordered spatial embedding to determine the ordered tree structure.

COROLLARY 5.5. *For $d \geq 8$, there exists a set $C^* \subseteq C(\mathbb{R}_+, \mathbb{R}_+) \times C(\mathbb{R}_+, \mathbb{R}^d)$ such that $\tilde{M}(C(\mathbb{R}_+, \mathbb{R}_+) \times C(\mathbb{R}_+, \mathbb{R}^d) \setminus C^*) = 0$ and, furthermore, if $(v, r), (\tilde v, \tilde r) \in C^*$, then*

$$r = \tilde r \quad \Leftrightarrow \quad (v, r) = (\tilde v, \tilde r).$$

PROOF. From the construction of $(v, r)$ and the definition of $\tilde{M}$ we can use Proposition 5.2 to deduce that there exists a set $C^*$ whose complement is $\tilde M$-null such that, if $(v, r) \in C^*$, the topological space $\mathcal{S} = r(\mathbb{R}_+)$ is a dendrite whose $c_d x^2 \ln\ln x^{-1}$-Hausdorff measure along arcs gives a metric $d_\mathcal{S}$ on $r(\mathbb{R}_+)$. The proof is completed on observing that we can also take as an assumption that on $C^*$ the function $v$ can be recovered via the relationship $v(t) = d_\mathcal{S}(0, r(t))$, for every $t \in \mathbb{R}_+$. □

**6. Brownian motion on spatial trees: quenched law.** Now that we have constructed the metric $d_\mathcal{S}$ on $\mathcal{S} \subseteq \mathbb{R}^d$ for $d \geq 8$, there is very little we have to do to build a canonical Markov process, $X^\mathcal{S}$, say, on $\mathcal{S}$ in high dimensions. We show in this section how the process $X^\mathcal{S}$, which we will call the Brownian motion on $\mathcal{S}$, can be obtained directly from $\mathcal{S}$ or, alternatively, it can be defined as $\phi(X^\mathcal{T})$, where $X^\mathcal{T}$ is a natural Markov process on the real tree $\mathcal{T}$. In the case when we normalize $\mu^\mathcal{T}$ (and $\mu^\mathcal{S}$) to have total mass one, we are able to describe these Markov processes as scaling limits of simple random walks on random graph trees embedded into Euclidean space; see Sections 8 and 9.

Let us start by introducing some known results about Dirichlet forms and Brownian motion on a compact real tree. Suppose $(\mathcal{T}, d_\mathcal{T})$ is a real tree and $\nu$ is a finite Borel measure on $\mathcal{T}$ that satisfies $\nu(A) > 0$ for every nonempty open set $A \subseteq \mathcal{T}$. Given a local, regular Dirichlet form $(\mathcal{E}_\mathcal{T}, \mathcal{F}_\mathcal{T})$ on $L^2(\mathcal{T}, \nu)$, we can use the standard association to define a nonnegative self-adjoint operator, $-\Delta_\mathcal{T}$, which has domain dense in $L^2(\mathcal{T}, \nu)$ and satisfies

$$\mathcal{E}_\mathcal{T}(f, g) = -\int_\mathcal{T} f \Delta_\mathcal{T} g \, d\nu \qquad \forall f \in \mathcal{F}_\mathcal{T}, g \in \mathcal{D}(\Delta_\mathcal{T}).$$

We can use this to define a reversible strong Markov process,

$$X^{\mathcal{T},\nu} = ((X_t^{\mathcal{T},\nu})_{t \geq 0}, \mathbf{P}_\sigma^{\mathcal{T},\nu}, \sigma \in \mathcal{T})$$

with semi-group given by $P_t := e^{t\Delta_\mathcal{T}}$. In fact, the locality of our Dirichlet form ensures that the process $X^{\mathcal{T},\nu}$ is a diffusion on $\mathcal{T}$. A fundamental example of a local, regular Dirichlet form is obtained as the electrical energy when we consider $(\mathcal{T}, d_\mathcal{T})$ to be an electrical network. In particular, the existence of a Dirichlet form for which the metric $d_\mathcal{T}$ describes the resistance between points of $\mathcal{T}$, so that

$$d_\mathcal{T}(x, y)^{-1} = \inf\{\mathcal{E}_\mathcal{T}(f, f) : f \in \mathcal{F}_\mathcal{T}, f(x) = 1, f(y) = 0\}$$

for every $x, y \in \mathcal{T}$, $x \neq y$, is guaranteed by Theorem 5.4 of [18]. The unique quadratic form with this property is known as the resistance form associated with $(\mathcal{T}, d_\mathcal{T})$ (see [19] for an introduction to resistance forms).

We follow Aldous [2] in defining a Brownian motion on $(\mathcal{T}, d_\mathcal{T}, \nu)$ to be a strong Markov process with continuous sample paths that is reversible with respect to its invariant measure $\nu$ and satisfies the following properties:

(i) For $\sigma_1, \sigma_2 \in \mathcal{T}$, $\sigma_1 \neq \sigma_2$, we have

$$\mathbf{P}_\sigma^{\mathcal{T},\nu}(T_{\sigma_1} < T_{\sigma_2}) = \frac{d_\mathcal{T}(b^\mathcal{T}(\sigma, \sigma_1, \sigma_2), \sigma_2)}{d_\mathcal{T}(\sigma_1, \sigma_2)} \qquad \forall \sigma \in \mathcal{T}$$

where $T_\sigma := \inf\{t > 0 : X_t^\mathcal{T} = \sigma\}$ is the hitting time of $\sigma \in \mathcal{T}$.

(ii) For $\sigma_1, \sigma_2 \in \mathcal{T}$, the mean occupation measure for the process started at $\sigma_1$ and killed on hitting $\sigma_2$ has density

$$2 d_\mathcal{T}(b^\mathcal{T}(\sigma, \sigma_1, \sigma_2), \sigma_2) \nu(d\sigma) \qquad \forall \sigma \in \mathcal{T}.$$

These properties guarantee the uniqueness of Brownian motion on $(\mathcal{T}, d_\mathcal{T}, \nu)$, and to construct the process we can use the resistance form described above. The following result can be proved using ideas from [9], Section 8.

PROPOSITION 6.1. *Let $(\mathcal{T}, d_\mathcal{T})$ be a compact real tree, $\nu$ be a finite Borel measure on $\mathcal{T}$ that satisfies $\nu(A) > 0$ for every nonempty open set $A \subseteq \mathcal{T}$, and $(\mathcal{E}_\mathcal{T}, \mathcal{F}_\mathcal{T})$ be the resistance form $(\mathcal{E}_\mathcal{T}, \mathcal{F}_\mathcal{T})$ associated with $(\mathcal{T}, d_\mathcal{T})$. Then $(\frac{1}{2}\mathcal{E}_\mathcal{T}, \mathcal{F}_\mathcal{T})$ is a local, regular Dirichlet form on $L^2(\mathcal{T}, \nu)$, and the corresponding Markov process $X^{\mathcal{T},\nu}$ is Brownian motion on $(\mathcal{T}, d_\mathcal{T}, \nu)$.*

From the above result, we are easily able to obtain the following. Note that, since we only consider one measure on $\mathcal{T}$ and one on $\mathcal{S}$, we henceforth drop the measure from the superscripts of the Brownian motions on these spaces.

PROPOSITION 6.2. *(a) For $\Theta$-a.e. $\mathcal{T}$, the Brownian motion $X^\mathcal{T}$ on the space $(\mathcal{T}, d_\mathcal{T}, \mu^\mathcal{T})$ exists.*
*(b) Let $d \geq 8$. For $M$-a.e. $(\mathcal{T}, \phi)$, the Brownian motion $X^\mathcal{S}$ on the space $(\mathcal{S}, d_\mathcal{S}, \mu^\mathcal{S})$ exists and, moreover, $X^\mathcal{S} = \phi(X^\mathcal{T})$.*

PROOF. In view of Propositions 5.2 and 6.1, it remains to prove that $X^\mathcal{S} = \phi(X^\mathcal{T})$, $M$-a.e. when $d \geq 8$. To do this, it will be sufficient to check that $\phi(X^\mathcal{T})$ satisfies the defining properties of Brownian motion on $(\mathcal{S}, d_\mathcal{S}, \mu^\mathcal{S})$. This is straightforward given the fact that $\phi$ is an isometry from $(\mathcal{T}, d_\mathcal{T})$ to $(\mathcal{S}, d_\mathcal{S})$ which satisfies $\mu^\mathcal{S} = \mu^\mathcal{T} \circ \phi^{-1}$. □

Since the map $\phi \colon \mathcal{T} \to \mathbb{R}^d$ is continuous for $M$-a.e. spatial tree $(\mathcal{T}, \phi)$, the law $\mathbf{P}_\rho^\mathcal{T} \circ \phi^{-1}$ of $\phi(X^\mathcal{T})$ is a well-defined probability measure on $C(\mathbb{R}_+, \mathbb{R}^d)$, $M$-a.e. Using the language of random walks in random environments, we say this is the quenched law of $\phi(X^\mathcal{T})$. Similarly, if $d \geq 8$, applying the fact that the identity map $I \colon (\mathcal{S}, d_\mathcal{S}) \to (\mathbb{R}^d, d_E)$ is continuous for $M$-a.e. spatial tree $(\mathcal{T}, \phi)$, we have that the law of the Brownian motion on $\mathcal{S}$, $\mathbf{P}_0^\mathcal{S}$, is a well-defined probability measure on $C(\mathbb{R}_+, \mathbb{R}^d)$, $M$-a.e., and we will call this the quenched law of the Brownian motion on $\mathcal{S}$. To allow us to define the annealed laws of $\phi(X^\mathcal{T})$ and $X^\mathcal{S}$ by averaging over the possible choices of environments, we need to show that there exists a probability space on which $\mathbf{P}_\rho^\mathcal{T} \circ \phi^{-1}$ and $\mathbf{P}_0^\mathcal{S}$ can be constructed measurably, and we will do this in the next section.

**7. Brownian motion on spatial trees: annealed law.** In the proofs of the measurability of the laws $\mathbf{P}_\rho^\mathcal{T} \circ \phi^{-1}$ and $\mathbf{P}_0^\mathcal{S}$, and the convergence results of later sections, it will be useful to approximate the spaces $\mathcal{T}$ and $\mathcal{S}$ by tree-like sets with only a finite number of branches. To this end, we introduce the concept of a (rooted ordered) graph spatial tree. This is a pair $(T, \phi)$, where $T$ is a (rooted) ordered finite graph tree with finite edge lengths, and $\phi$ is a continuous $\mathbb{R}^d$-valued map whose domain is the real tree $\overline{T}$ naturally associated with $T$ by adding line segments to $T$ along its edges and extending the graph distance on $T$ to a metric $d_{\overline{T}}$ on $\overline{T}$ in the natural way so that the line segment corresponding to an edge with weight $|e|$ is isometric to $[0, |e|]$. We will assume that $\phi$ maps the root of $T$ to the origin in $\mathbb{R}^d$. Note that for each graph spatial tree the pair $(\overline{T}, \phi)$ is an element of $\mathbb{T}_{\rm sp}$, and, since $T$ is finite, $\overline{T}$ is compact. Moreover, the fact that $T$ is finite means that we can define a probability measure $\lambda^{\overline{T}}$ on $\overline{T}$ to be the renormalized Lebesgue measure (so that the $\lambda^{\overline{T}}$-measure of a line segment in $\overline{T}$ is proportional to its length), and by Proposition 6.1, there is no problem in defining the Brownian motion on $(\overline{T}, \lambda^{\overline{T}})$.

The topology we consider on the space of graph spatial trees is a generalization of the topology considered in Section 4 of [8]. For $(T, \phi)$ a graph spatial tree, write $T = (T^*; |e_1|, \ldots, |e_l|)$, where $T^*$ represents the "shape" of $T$ (the ordered graph tree without edge lengths) and $|e_1|, \ldots, |e_l|$ represents the collection of edge lengths. Then, if $(T, \phi)$ and $(T', \phi')$ are two graph spatial trees, define a distance $d_1$ between $T$ and $T'$ by setting $d_1(T, T') := \infty$, when $T^* \neq T'^*$, and

$$d_1(T, T') := \sup_i ||e_i| - |e_i'||$$

otherwise. When $T^* = T'^*$, we have a homeomorphism $\Upsilon_{\overline{T}, \overline{T'}} : \overline{T} \to \overline{T'}$, under which the point $x \in \overline{T}$, which is a distance $\alpha$ along an edge $e$ (considered from the vertex at the end of $e$ which is closest to the root), is mapped to the point $x' \in \overline{T'}$ which is a distance $|e'|\alpha/|e|$ along the corresponding edge $e'$, and so we can define

$$d_2(\phi, \phi') := \sup_{x \in \overline{T}} d_E(\phi(x), \phi'(\Upsilon_{\overline{T}, \overline{T'}}(x)))$$

which yields a metric $d_0((T, \phi), (T', \phi')) := (d_1(T, T') + d_2(\phi, \phi')) \wedge 1$. It is straightforward to check that, when equipped with the topology induced by this metric, the collection of graph spatial trees is separable. In [8], it was shown that if $(T_n)_{n \geq 1}$ is a sequence of ordered graph trees that converge with respect to the distance $d_1$ to an ordered graph tree $T$, then $\Upsilon_{\overline{T}_n, \overline{T}}(B^n)$, where $B^n = (B^n_t)_{t \geq 0}$ is the Brownian motion on $(\overline{T}_n, \lambda^{\overline{T}_n})$ started from the root, converges in distribution in the space $C(\mathbb{R}_+, \overline{T})$ to $B$, the Brownian motion on $(\overline{T}, \lambda^{\overline{T}})$ started from the root. By mapping this result into $\mathbb{R}^d$ in the obvious way, we are easily able to deduce from this the following lemma.

LEMMA 7.1. *Suppose that $\{(T_n, \phi_n)\}_{n \geq 1}$ is a sequence of graph spatial trees that converge with respect to the metric $d_0$ to a graph spatial tree $(T, \phi)$. For each $n$, let $B^n = (B^n_t)_{t \geq 0}$ be the Brownian motion on $(\overline{T}_n, \lambda^{\overline{T}_n})$ started from the root, and let $B = (B_t)_{t \geq 0}$ be the Brownian motion on $(\overline{T}, \lambda^{\overline{T}})$ started from the root, then $(\phi_n(B^n))_{n \geq 1}$ converges in distribution to $\phi(B)$ in the space $C(\mathbb{R}_+, \mathbb{R}^d)$.*

Let us continue by considering a vector $(\sigma_1, \ldots, \sigma_k)$ of elements of a real tree $\mathcal{T}$. Define the reduced subtree $\mathcal{T}(\sigma_1, \ldots, \sigma_k)$ to be the graph tree with vertex set

$$V(\mathcal{T}(\sigma_1, \ldots, \sigma_k)) := \{b^{\mathcal{T}}(\sigma, \sigma', \sigma'') : \sigma, \sigma', \sigma'' \in \{\rho, \sigma_1, \ldots, \sigma_k\}\},$$

and graph tree structure induced by the arcs of $\mathcal{T}$, so that two elements $\sigma$ and $\sigma'$ of $V(\mathcal{T}(\sigma_1,\ldots,\sigma_k))$ are connected by an edge if and only if $\sigma \neq \sigma'$ and also $\Gamma^{\mathcal{T}}_{\sigma,\sigma'} \cap V(\mathcal{T}(\sigma_1,\ldots,\sigma_k)) = \{\sigma,\sigma'\}$. We set the length of an edge $\{\sigma,\sigma'\}$ to be equal to $d_{\mathcal{T}}(\sigma,\sigma')$. Furthermore, we can use the order of the vector $(\sigma_1,\ldots,\sigma_k)$ to induce an ordering of vertices in the graph tree $\mathcal{T}(\sigma_1,\ldots,\sigma_k)$. If $(\mathcal{T},\phi) \in \mathbb{T}_{\mathrm{sp}}$, it is possible to restrict $\phi$ to $\overline{\mathcal{T}(\sigma_1,\ldots,\sigma_k)}$ to obtain a graph spatial tree $(\mathcal{T}(\sigma_1,\ldots,\sigma_k),\phi)$. In the following result, by considering the Brownian motions on an increasing sequence of reduced subtrees of graph spatial trees, we show that the law $\mathbf{P}^{\mathcal{T}}_{\rho} \circ \phi^{-1}$ is a measurable function of the tour defining $(\mathcal{T},\phi)$.

PROPOSITION 7.2. *The map from the tour $(v,r)$ to $\mathbf{P}^{\mathcal{T}}_{\rho} \circ \phi^{-1}$ is measurable with respect to the $\tilde{M}$-completion of the standard topology on the space $C(\mathbb{R}_+,\mathbb{R}_+) \times C(\mathbb{R}_+,\mathbb{R}^d)$ and the topology induced by the weak convergence of probability measures on $C(\mathbb{R}_+,\mathbb{R}^d)$.*

PROOF. Let $\{(v_n,r_n,u_n)\}_{n\geq 1}$ be a sequence in $C(\mathbb{R}_+,\mathbb{R}_+) \times C(\mathbb{R}_+,\mathbb{R}^d) \times [0,1]^{\mathbb{N}}$ converging to $(v,r,u)$, where $\{(v_n,r_n)\}_{n\geq 1}$ and $(v,r)$ are tours. We will write $u_n = (u_n^m)_{m\geq 1}$ and $u = (u^m)_{m\geq 1}$, and rescaled versions of these by $\tilde{u}_n = (\tau(v_n)u_n^m)_{m\geq 1}$ and $\tilde{u} = (\tau(v)u^m)_{m\geq 1}$, where $\tau(\cdot)$ is the length of excursion function, as defined in Section 2.1. By [24], Theorem 2.1, we also have that the corresponding sequence of snakes $\{(v_n,w_n)\}_{n\geq 1}$ converges to a limit snake, $(v,w)$, say. This implies, for every fixed $k \geq 1$, that the vector

$$(v_n(\tilde{u}_n^{(1)}), w_n(\tilde{u}_n^{(1)}), m_{v_n}(\tilde{u}_n^{(1)},\tilde{u}_n^{(2)}),\ldots,$$
$$v_n(\tilde{u}_n^{(2)}), w_n(\tilde{u}_n^{(2)}),\ldots,m_{v_n}(\tilde{u}_n^{(k-1)},\tilde{u}_n^{(k)}), v_n(\tilde{u}_n^{(k)}), w_n(\tilde{u}_n^{(k)})),$$

where $(\tilde{u}_n^{(m)})_{m=1}^k$ is a nondecreasing ordering of $(\tilde{u}_n^m)_{m=1}^k$, converges to

$$(v(\tilde{u}^{(1)}), w(\tilde{u}^{(1)}), m_v(\tilde{u}^{(1)},\tilde{u}^{(2)}), v(\tilde{u}^{(2)}), w(\tilde{u}^{(2)}),\ldots,$$
$$m_v(\tilde{u}^{(k-1)},\tilde{u}^{(k)}), v(\tilde{u}^{(k)}), w(\tilde{u}^{(k)})),$$

where $(\tilde{u}^{(m)})_{m=1}^k$ is a nondecreasing ordering of $(\tilde{u}^m)_{m=1}^k$. If we set

(7.1) $$T^{(k)}_{v,u} := \mathcal{T}_v([\tilde{u}^{(1)}],\ldots,[\tilde{u}^{(k)}]),$$

where $[t]$ represents the equivalence classes under the equivalence defined at (2.5), and define $T^{(k)}_{v_n,u_n}$ similarly, then it follows that, if $T^{(k)}_{v,u}$ has no vertex of degree greater than three, then the graph spatial tree $(T^{(k)}_{v_n,u_n},\phi_{v_n,r_n})$ converges to $(T^{(k)}_{v,u},\phi_{v,r})$ with respect to the metric $d$ (cf. the proof of [3], Theorem 20); consequently, by Lemma 7.1, the law of $\phi_{v_n,r_n}(B^{n,k})$, where $B^{n,k}$ is the Brownian motion on $(\overline{T}^{(k)}_{v_n,u_n}, \lambda^{\overline{T}^{(k)}_{v_n,u_n}})$ started from the root, converges to the law, $\mathbf{P}_0^{\mathcal{S}(k)}$, say, of $\phi_{v,r}(B^{(k)})$, where $B^{(k)}$ is the Brownian motion on $(\overline{T}^{(k)}_{v,u}, \lambda^{\overline{T}^{(k)}_{v,u}})$ started from the root.

It is known (see [14], Theorem 4.6) that, for $N$-a.e. realization of $v$, the set $\mathcal{T}_v \setminus \{x\}$ has at most three connected components for any $x \in \mathcal{T}_v$. Hence, by applying the conclusion of the previous paragraph, we are able to deduce that there exists a measurable set $\Gamma \subseteq C(\mathbb{R}_+, \mathbb{R}_+) \times C(\mathbb{R}_+, \mathbb{R}^d) \times [0,1]^\mathbb{N}$ with $\tilde{M} \otimes \lambda^{\otimes \mathbb{N}}_{[0,1]}(\Gamma^c) = 0$, where $\lambda_{[0,1]}$ is the Lebesgue measure on $[0,1]$, such that the map from $(v,r,u) \in \Gamma$ to $\mathbf{P}_0^{\mathcal{S}(k)}$ is continuous on $\Gamma$, and therefore measurable on $C(\mathbb{R}_+, \mathbb{R}_+) \times C(\mathbb{R}_+, \mathbb{R}^d) \times [0,1]^\mathbb{N}$ with respect to the $\tilde{M} \otimes \lambda^{\otimes \mathbb{N}}_{[0,1]}$-completion of the standard product topology on this space.

By following the proof of [8], Lemma 3.1, we obtain that, for $\tilde{M} \otimes \lambda^{\otimes \mathbb{N}}_{[0,1]}$-a.e. realization of $(v,r,u)$, the law of $B^{(k)}$ converges weakly in the space of Borel probability measures on $C(\mathbb{R}_+, \mathcal{T}_v)$ to the law $\mathbf{P}^{\mathcal{T}_v}_{\rho_v}$. Hence, by the continuity of $\phi_{v,r}$, we have that $\mathbf{P}_0^{\mathcal{S}(k)}$ converges weakly to $\mathbf{P}^{\mathcal{T}_v}_{\rho_v} \circ \phi^{-1}_{v,r}$, $\tilde{M} \otimes \lambda^{\otimes \mathbb{N}}_{[0,1]}$-a.e. Since a limit of measurable functions is again measurable, the map from $(v,r,u)$ to $\mathbf{P}^{\mathcal{T}_v}_{\rho_v} \circ \phi^{-1}_{v,r}$ is measurable with respect to the topology described in the previous paragraph. Noting that

$$\mathbf{P}^{\mathcal{T}_v}_{\rho_v} \circ \phi^{-1}_{v,r} = \int_{[0,1]^\mathbb{N}} \mathbf{P}^{\mathcal{T}_v}_{\rho_v} \circ \phi^{-1}_{v,r} \lambda^{\otimes \mathbb{N}}_{[0,1]}(du),$$

the result follows. $\square$

As an immediate consequence of the above result, it is possible to define a measure $\mathbb{M}$ on $C(\mathbb{R}_+, \mathbb{R}_+) \times C(\mathbb{R}_+, \mathbb{R}^d)^2$ which satisfies

(7.2) $$\mathbb{M}(A \times B) := \int_A \mathbf{P}^{\mathcal{T}}_\rho(\phi(X^{\mathcal{T}}) \in B) \tilde{M}(d(v,r)),$$

for measurable $A \subseteq C(\mathbb{R}_+, \mathbb{R}_+) \times C(\mathbb{R}_+, \mathbb{R}^d)$ and $B \subseteq C(\mathbb{R}_+, \mathbb{R}^d)$. This represents first choosing a tour $(v,r)$ by the measure $\tilde{M}$ [which means that the resulting spatial tree $(\mathcal{T}, \phi)$ has marginal $M$], and then observing the Brownian motion on the real tree $\mathcal{T}$ mapped into Euclidean space by $\phi$; in the random walk in random environment terminology, the law of $\phi(X^{\mathcal{T}})$ under $\mathbb{M}$ is the annealed law of $\phi(X^{\mathcal{T}})$. For $d \geq 8$, we can apply Proposition 6.2 to simplify the formula at (7.2) so that the integrand only depends on the set $\mathcal{S}$ rather than the whole spatial tree $(\mathcal{T}, \phi)$. In particular, $\mathbb{M}$ satisfies

$$\mathbb{M}(A \times B) := \int_A \mathbf{P}^{\mathcal{S}}_0(X^{\mathcal{S}} \in B) \tilde{M}(d(v,r))$$

for measurable $A \subseteq C(\mathbb{R}_+, \mathbb{R}_+) \times C(\mathbb{R}_+, \mathbb{R}^d)$ and $B \subseteq C(\mathbb{R}_+, \mathbb{R}^d)$. In this high-dimensional case, we call the law of $X^{\mathcal{S}}$ under $\mathbb{M}$ the annealed law of the Brownian motion on $\mathcal{S}$.

**8. Quenched convergence.** The aim of this section is to prove convergence results for the simple random walks on a family of graph spatial trees, given that the associated discrete tours converge to a typical realization of the normalized Brownian tour. In particular, we consider a family $\{(T_n, \phi_n)\}_{n\geq 1}$ of graph spatial trees, as defined in the previous section, such that each graph $T_n$ has $n$ vertices and is unweighted, by which we mean that each edge has length one. We define $\mu^{T_n}$ to be the uniform probability measure on the vertices of $T_n$ and, analogous to the definitions of $\mathcal{S}$ and $\mu^{\tilde{\mathcal{S}}}$, set

$$S_n := \phi_n(\overline{T}_n), \qquad \mu^{S_n} := \mu^{T_n} \circ \phi_n^{-1}.$$

Furthermore, let $X^{T_n} = ((X_m^{T_n})_{m\geq 0}, \mathbf{P}_x^{T_n}, x \in T_n)$ be the usual discrete time simple random walk on the vertices of $T_n$. To define $X^{T_n}$ at all positive times, we linearly interpolate between integers (for this to make sense, we suppose that the walk takes values in the real tree version $\overline{T}_n$ of $T_n$ obtained by adding unit line segments along edges).

Since the excursion description of ordered graph trees and the corresponding discrete tour and snake description are well documented in [3] and [24] respectively, we will not present the full details, but simply highlight the results that will be important here. Define $\tilde{V}_n : \{1, \ldots, 2n-1\} \to T_n$ to be the depth-first search around the vertices of the ordered graph tree $T_n$, starting from the root at time one. Extend this map to the interval $[0, 2n]$ by setting $\tilde{V}_n(0) = \tilde{V}_n(2n) = \rho_n$, where $\rho_n$ is the root of $T_n$, and linearly interpolating (similarly to the extension of the simple random walk, we now consider that $\tilde{V}_n$ takes values in the real tree $\overline{T}_n$). The search-depth function $V_n \in C([0,1], \mathbb{R}_+)$ is given by, for $t \in [0,1]$,

$$V_n(t) := d_{\overline{T}_n}(\rho_n, \tilde{V}_n(2nt)),$$

where $d_{\overline{T}_n}$ is the metric on $\overline{T}_n$. A related function in $\mathbb{R}^d$ is given by, for $t \in [0,1]$,

$$R_n(t) := \phi_n(\tilde{V}_n(2nt)),$$

which is the discrete head process and is an element of $C([0,1], \mathbb{R}^d)$. The process $(V_n, R_n)$ is the discrete tour associated with $(T_n, \phi_n)$, although from now on we will commonly refer to the normalized discrete tour, which is defined by setting

$$(v_n, r_n) := (n^{-1/2} V_n, n^{-1/4} R_n).$$

It is easy to check that the normalized discrete tour $(v_n, r_n)$ contains all the information about the graph spatial tree $(T_n, \phi_n)$. The corresponding normalized discrete snake $w_n$ is a continuous function taking its values in the space of $\mathbb{R}^d$-valued stopped paths, and is defined to satisfy

$$w_n(t)(s) := n^{-1/4} \phi_n(\gamma_{\rho_n, \tilde{V}_n(2nt)}^{\overline{T}_n}(n^{1/2} s))$$

for $s \leq v_n(t)$, and $w_n(t)(s) = r_n(t)$ otherwise.

We can now state and prove the main result of this section. Note that, for $(\mathcal{T}, \phi) \in \mathbb{T}_{\mathrm{sp}}$, we write $(\alpha \mathcal{T}, \beta \phi)$ to represent the real tree $(\mathcal{T}, \alpha d_{\mathcal{T}}, \rho)$ and map $\sigma \mapsto \beta \phi(\sigma)$, for $\sigma \in \mathcal{T}$; graph spatial trees will be rescaled similarly. The definition of the measure $\tilde{M}^{(1)}$ should be recalled from Section 2.3.

THEOREM 8.1. *There exists a set $C^* \subseteq C([0,1], \mathbb{R}_+) \times C([0,1], \mathbb{R}^d)$ with $\tilde{M}^{(1)}(C^*) = 1$ such that, if $(v_n, r_n) \to (v, r)$ in $C([0,1], \mathbb{R}_+) \times C([0,1], \mathbb{R}^d)$ for some $(v, r) \in C^*$, then the following statements hold, where $(\mathcal{T}, \phi)$ is the spatial tree associated with $(v, r)$:*

(a) $(n^{-1/2} \overline{T}_n, n^{-1/4} \phi_n) \to (\mathcal{T}, \phi)$ *in the space* $\mathbb{T}_{\mathrm{sp}}$.

(b) $n^{-1/4} S_n \to \mathcal{S}$ *with respect to the Hausdorff topology on compact subsets of* $\mathbb{R}^d$.

(c) $\mu^{S_n}(n^{1/4} \cdot) \to \mu^{\mathcal{S}}$ *weakly as Borel probability measures on* $\mathbb{R}^d$.

(d) $(n^{-1/4} \phi_n(X^{T_n}_{tn^{3/2}}))_{t \geq 0} \to \phi(X^{\mathcal{T}})$ *in distribution in* $C(\mathbb{R}_+, \mathbb{R}^d)$.

PROOF. Assume that $(v_n, r_n) \to (v, r)$ in $C([0,1], \mathbb{R}_+) \times C([0,1], \mathbb{R}^d)$. For each $n$, define a correspondence $\mathcal{C}_n$ between $n^{-1/2} \overline{T}_n$ and $\mathcal{T}$ by

$$\mathcal{C}_n := \{(\sigma, \sigma') : \sigma = \tilde{V}_n(2nt), \sigma' = [t], \text{ for some } t \in [0,1]\},$$

where $[t]$ represents the equivalence classes of $[0,1]$ under the equivalence defined at (2.5), and, since $n^{-1/2} \overline{T}_n := (\overline{T}_n, n^{-1/2} d_{\overline{T}_n}, \rho_n)$, we note that the function $\tilde{V}_n$ can indeed be considered as a function from $[0, 2n]$ to $n^{-1/2} \overline{T}_n$. Similarly to the proof of Proposition 2.4, this correspondence allows us to deduce that

$$d_{\mathrm{sp}}((n^{-1/2} \overline{T}_n, n^{-1/4} \phi_n), (\mathcal{T}, \phi)) \leq 4 \|v_n - v\|_\infty + \|r_n - r\|_\infty$$

and, therefore, part (a) holds. As noted in Section 2.2, the map $(\mathcal{T}, \phi) \mapsto \phi(\mathcal{T})$ is continuous, hence, part (b) is an immediate consequence of part (a). To prove part (c), we start by considering the Lebesgue measure $\lambda_{[0,1]}$ on $[0,1]$. By the characterization of $\mu^{\mathcal{T}}$ at (2.6), it is clear that $\lambda_{[0,1]} \circ r^{-1}$ is identical to $\mu^{\mathcal{T}} \circ \phi^{-1} = \mu^{\mathcal{S}}$. For graph trees, the analogous representation is not quite as straightforward, because the uniform measure on $[0,1]$ does not map to the uniform measure on the vertices on $T_n$ in such a simple way. However, this problem is not major. Define the function $\alpha_n : [0,1] \to [0,1]$ by setting

$$\alpha_n(t) := \begin{cases} \lfloor 2nt \rfloor / 2n, & \text{if } v_n(\lfloor 2nt \rfloor / 2n) \geq v_n(\lceil 2nt \rceil / 2n), \\ \lceil 2nt \rceil / 2n, & \text{otherwise.} \end{cases}$$

It is clear from the definition that $\sup_{t \in [0,1]} |t - \alpha_n(t)| \leq 1/2n$, regardless of the value of $v_n$. Furthermore, by applying an argument similar to Lemma 12 of [3], it is possible to show that if $U$ is a random variable with law $\lambda_{[0,1]}$, then $\tilde{V}_n(2n\alpha_n(U))$ is uniform on the vertices of $T_n$. Since by assumption $n^{-1/4}\phi_n(\tilde{V}_n(2n\alpha_n(U))) = r_n(\alpha_n(U)) \to r(U)$, it follows that $\mu^{S_n}(n^{1/4}\cdot) \to \mu^{\mathcal{S}}$, as required for part (c) to hold.

To prove part (d), we will use the idea of reduced subtrees, as in the proof of Proposition 7.2. First, note that [24], Theorem 2.1, implies that $(v_n, w_n)$ converges to $(v, w)$, where $w_n$ is the normalized discrete snake associated with $(v_n, r_n)$ and $w$ is the snake associated with $(v, r)$. Thus, we can proceed similarly to the proof of Proposition 7.2 to deduce that if we let $u = (u^m)_{m\geq 1}$ be a sequence taking values in $[0, 1]$, define $T^{(k)} = T_{v,u}^{(k)}$ as at (7.1), and introduce a reduced subtree of $\overline{T}_n$ by setting

$$T_n^{(k)} := \overline{T}_n(\tilde{V}_n(2n\alpha_n(u^{(1)})), \ldots, \tilde{V}_n(2n\alpha_n(u^{(k)}))),$$

where $(u^{(m)})_{m=1}^k$ is a nondecreasing ordering of $(u^m)_{m=1}^k$, then

(8.1) $$(n^{-1/2}T_n^{(k)}, n^{-1/4}\phi_n) \to (T^{(k)}, \phi)$$

in the space of graph spatial trees whenever $T^{(k)}$ has no vertex of degree greater than three. Recalling that (see [14], Theorem 4.6), for $N$-a.e. realization of $v$, the set $\mathcal{T}_v \setminus \{x\}$ has at most three connected components for any $x \in \mathcal{T}_v$, we can take the convergence of the previous sentence as an assumption. Consequently, if we let

$$X^{T_n^{(k)}} = (X_t^{T_n^{(k)}})_{t\geq 0}$$

be the nearest neighbor discrete time simple random walk on the vertices of $T_n$ contained in $\overline{T}_n^{(k)}$, extended to a continuous time process taking values in $\overline{T}_n^{(k)}$ by linear interpolation, then [8], Lemma 4.2, implies that, for $N \otimes \lambda_{[0,1]}^{\otimes \mathbb{N}}$-a.e. realization of $(v, u)$,

$$\Upsilon_{n^{-1/2}\overline{T}_n^{(k)}, \overline{T}^{(k)}}(X_{tn\Lambda_n^{(k)}}^{T_n^{(k)}}) \to B^{(k)}$$

as $n \to \infty$, in distribution in $C(\mathbb{R}_+, \overline{T}^{(k)})$, where $\Lambda_n^{(k)}$ is the sum of the lengths of the edges of $T_n^{(k)}$, and $B^{(k)}$ is defined, as in the proof of Proposition 7.2, to be the Brownian motion on $(\overline{T}^{(k)}, \lambda^{\overline{T}^{(k)}})$ started from the root. Applying (8.1), it is possible to deduce by mapping this conclusion into $\mathbb{R}^d$ using $\phi$ that, for $\tilde{M}^{(1)} \otimes \lambda_{[0,1]}^{\otimes \mathbb{N}}$-a.e. realization of $(v, r, u)$,

(8.2) $$n^{-1/4}\phi_n(X_{tn\Lambda_n^{(k)}}^{T_n^{(k)}}) \to \phi(B^{(k)})$$

in distribution in $C(\mathbb{R}_+, \mathbb{R}^d)$. Note that, as in the proof of Proposition 7.2, we also have that

(8.3) $$\phi(B^{(k)}) \to \phi(X^{\mathcal{T}})$$

as $k \to \infty$, in distribution in $C(\mathbb{R}_+, \mathbb{R}^d)$, for $\tilde{M}^{(1)} \otimes \lambda_{[0,1]}^{\otimes \mathbb{N}}$-a.e. realization of $(v, r, u)$.

Finally, (cf. [8], Proposition 7.1) we can define $X^{T_n}$ and $X^{T_n^{(k)}}$ on the same probability space (with probability measure $\mathbf{P}$) in such a way that

$$(8.4) \quad \lim_{k \to \infty} \limsup_{n \to \infty} \mathbf{P}\left( n^{-1/2} \sup_{t \in [0, t_0]} d_{\overline{T}_n}(X^{T_n}_{tn^{3/2}}, X^{T_n^{(k)}}_{tn\Lambda_n^{(k)}}) > \varepsilon \right) = 0,$$

for every $\varepsilon, t_0 > 0$. The desired conclusion will follow easily from (8.2) and (8.3) by applying [5], Theorem 3.2, if we are able to replace the above tightness condition by

$$(8.5) \quad \lim_{k \to \infty} \limsup_{n \to \infty} \mathbf{P}\left( n^{-1/4} \sup_{t \in [0, t_0]} d_E(\phi_n(X^{T_n}_{tn^{3/2}}), \phi_n(X^{T_n^{(k)}}_{tn\Lambda_n^{(k)}})) > \varepsilon \right) = 0,$$

for every $\varepsilon, t_0 > 0$. First, for $\delta > 0$, we have

$$\sup_{d_{\overline{T}_n}(\sigma_1, \sigma_2) < \delta n^{1/2}} n^{-1/4} d_E(\phi_n(\sigma_1), \phi_n(\sigma_2)) = \sup_{d_{v_n}(s, t) < \delta} d_E(r_n(s), r_n(t)),$$

where $d_{v_n}$ is defined by the formula at (2.4). By assumption, this expression converges to

$$\sup_{d_v(s,t) < \delta} d_E(r(s), r(t)) = \sup_{d_{\mathcal{T}}(\sigma_1, \sigma_2) < \delta} d_E(\phi(\sigma_1), \phi(\sigma_2))$$

as $n \to \infty$. Since $\phi$ is continuous on $\mathcal{T}$, which is compact, it follows that

$$(8.6) \quad \lim_{\delta \to 0} \limsup_{n \to \infty} \sup_{d_{\overline{T}_n}(\sigma_1, \sigma_2) < \delta n^{1/2}} n^{-1/4} d_E(\phi_n(\sigma_1), \phi_n(\sigma_2)) = 0.$$

Combining this with (8.4), we obtain (8.5). □

As the subsequent theorem demonstrates, the convergence of reduced subtrees deduced in the above proof is intrinsically linked with the convergence of tours. This result is a quenched generalization of [3], Theorem 20, which details conditions for the convergence of the search-depth processes of random ordered graph trees to the normalized Brownian excursion; we will prove the corresponding annealed version in the next section. To state our result, we continue to use the notation $T_n^{(k)}$ and $T^{(k)}$ introduced in the proof of Theorem 8.1. For a compact subset $A$ of a compact real tree $\mathcal{T}$, we set

$$\Delta(\mathcal{T}, A) := \sup_{\sigma_1 \in \mathcal{T}} \inf_{\sigma_2 \in A} d_{\mathcal{T}}(\sigma_1, \sigma_2),$$

which measures the usual Hausdorff distance between $A$ and $\mathcal{T}$. Finally, we define $\phi_n^{(k)}$ to be the map from $\overline{T}_n^{(k)}$ to $\mathbb{R}^d$ which is equal to $\phi_n$ on the graph vertices of $T_n^{(k)}$ and linear along the line-segments between them, and define $\phi^{(k)} : \overline{T}^{(k)} \to \mathbb{R}^d$ similarly.

THEOREM 8.2. *There exists a set* $D^* \subseteq C([0,1], \mathbb{R}_+) \times C([0,1], \mathbb{R}^d) \times [0,1]^{\mathbb{N}}$ *with* $\tilde{M}^{(1)} \otimes \lambda_{[0,1]}^{\otimes \mathbb{N}}(D^*) = 1$ *such that if* $(v, r, u) \in D^*$, *then the following three conditions are equivalent:*

(a) $(v_n, r_n) \to (v, r)$ *in the space* $C([0,1], \mathbb{R}_+) \times C([0,1], \mathbb{R}^d)$.

(b) *The convergence at (8.1) holds for each* $k \in \mathbb{N}$. *Furthermore,*

$$(8.7) \qquad \lim_{k \to \infty} \limsup_{n \to \infty} n^{-1/2} \Delta(\overline{T}_n, \overline{T}_n^{(k)}) = 0$$

*and also (8.6) is satisfied.*

(c) *The convergence at (8.1) holds for each* $k \in \mathbb{N}$ *when* $\phi_n$, $\phi$ *are replaced by* $\phi_n^{(k)}$, $\phi^{(k)}$ *respectively. Furthermore, (8.6) and (8.7) are satisfied.*

PROOF. The existence of a set $D^*$ with $\tilde{M}^{(1)} \otimes \lambda_{[0,1]}^{\otimes \mathbb{N}}(D^*) = 1$ upon which condition (a) implies (8.1) and (8.6) was demonstrated in the proof of the previous result. To prove (8.7), we can apply a deterministic version of the proof of [3], Theorem 20. Thus, (a) implies (b), from which (c) follows easily.

To prove (a) from (c), we start by noting that the convergence of subtrees implies that

$$(v_n(u^{(1)}), r_n(u^{(1)}), \ldots, v_n(u^{(k)}), r_n(u^{(k)}))$$
$$\to (v(u^{(1)}), r(u^{(1)}), \ldots, v(u^{(k)}), r(u^{(k)})),$$

where $(u^{(m)})_{m=1}^k$ is a nondecreasing ordering of $(u^m)_{m=1}^k$. Since we can take as an assumption that on $D^*$ the sequence $u$ is dense in $[0,1]$, to complete the proof it remains to show that $(v_n)_{n \geq 1}$ is tight in $C([0,1], \mathbb{R}_+)$ and $(r_n)_{n \geq 1}$ is tight in $C([0,1], \mathbb{R}^d)$ whenever (c) holds. To obtain the tightness of $(v_n)_{n \geq 1}$, we can again apply a deterministic version of the proof of [3], Theorem 20. Finally, we note that

$$\sup_{|s-t|<\delta} d_E(r_n(s), r_n(t)) \leq \sup_{d_{\overline{T}_n}(\sigma_1, \sigma_2) < 3\varepsilon(n,\delta) n^{1/2}} n^{-1/4} d_E(\phi_n(\sigma_1), \phi_n(\sigma_2)),$$

where $\varepsilon(n, \delta) := \sup_{|s-t|<\delta} |v_n(s) - v_n(t)|$. Applying this bound, the tightness of $(v_n)_{n \geq 1}$ and (8.6), it is an elementary exercise to deduce the tightness of $(r_n)_{n \geq 1}$. $\square$

To complete this section, let us briefly comment on the difference between (b) and (c) in the above theorem. First, observe that $\phi_n$ could be any continuous function on the edges of $T_n^{(k)}$, whereas $\phi_n^{(k)}$ simply records the increments of $\phi_n$ along the edges. Thus, condition (b) requires that the image under $\phi_n$ of an edge converges to a (typical) segment of a Brownian motion path in $\mathbb{R}^d$. In contrast, condition (c) requires the weaker condition that the increment of $\phi_n$ along each edge converges to the corresponding Brownian increment.

**9. Annealed convergence.** Applying the measurability of the map $(v, r) \to \mathbf{P}_\rho^\mathcal{T} \circ \phi^{-1}$ and the quenched convergence result that was proved in the previous section, we are able to establish a distributional convergence property for the simple random walks on a sequence of random graph spatial trees whose normalized discrete tours converge in distribution to the normalized Brownian tour. More specifically, we start by assuming that for each $n \in \mathbb{N}$ we have a probability measure $\tilde{M}_n$ on normalized discrete tours such that if $(v_n, r_n)$ is in the support of $\tilde{M}_n$, then the graph tree $T_n$ corresponding to $v_n$ has $n$ vertices and is unweighted. We can subsequently define a probability measure $\mathbb{M}_n$ on $C([0,1], \mathbb{R}_+) \times C([0,1], \mathbb{R}^d) \times C(\mathbb{R}_+, \mathbb{R}^d)$ that satisfies

$$(9.1) \quad \mathbb{M}_n(A \times B) = \int_A \mathbf{P}_{\rho_n}^{T_n}((n^{-1/4}\phi_n(X_{tn^{3/2}}^{T_n}))_{t \geq 0} \in B)\tilde{M}_n(d(v_n, r_n)),$$

for measurable $A \subseteq C([0,1], \mathbb{R}_+) \times C([0,1], \mathbb{R}^d)$ and $B \subseteq C(\mathbb{R}_+, \mathbb{R}^d)$ (the necessary measurability of the simple random walk laws is easily checked). Given these measures, which are the annealed measures of the normalized discrete tours and the associated simple random walks embedded into $\mathbb{R}^d$, the main result of this section is that if the laws of the discrete tours $\tilde{M}_n$ converge to the law of the normalized Brownian tour $\tilde{M}^{(1)}$, then $\mathbb{M}_n$ converges to $\mathbb{M}^{(1)}$, where we assume that $\mathbb{M}^{(1)}$ is a probability measure on $C([0,1], \mathbb{R}_+) \times C([0,1], \mathbb{R}^d) \times C(\mathbb{R}_+, \mathbb{R}^d)$ defined similarly to the annealed law $\mathbb{M}$ with $\tilde{M}$ replaced by $\tilde{M}^{(1)}$ in (7.2) (to justify this replacement, we note that it is straightforward to check that Proposition 7.2 holds when $\tilde{M}$ is replaced by $\tilde{M}^{(1)}$).

THEOREM 9.1. *If $\tilde{M}_n \to \tilde{M}^{(1)}$ weakly as probability measures on the space $C([0,1], \mathbb{R}_+) \times C([0,1], \mathbb{R}^d)$, then $\mathbb{M}_n \to \mathbb{M}^{(1)}$ weakly as probability measures on the space $C([0,1], \mathbb{R}_+) \times C([0,1], \mathbb{R}^d) \times C(\mathbb{R}_+, \mathbb{R}^d)$.*

PROOF. Following the proof of [8], Theorem 1.2, it is elementary to check that this result is a consequence of Theorem 8.1. □

Before continuing, we remark that Theorem 8.1 allows us to deduce that the convergence of $\tilde{M}_n \to \tilde{M}^{(1)}$ also implies the convergence of the laws of the sets $n^{-1/4}S_n$ and measures $\mu^{S_n}(n^{1/4}\cdot)$ under $\tilde{M}_n$ to the laws of $\mathcal{S}$ and $\mu^\mathcal{S}$, respectively, under $\tilde{M}^{(1)}$. To complete this section, we state the annealed version of Theorem 8.2, which can be proved by making the obvious changes to the proof of Theorem 8.2.

THEOREM 9.2. *The following three conditions are equivalent:*

(a) *$\tilde{M}_n \to \tilde{M}^{(1)}$ weakly as probability measures on the space $C([0,1], \mathbb{R}_+) \times C([0,1], \mathbb{R}^d)$.*

(b) *For each $k \in \mathbb{N}$, then*

$$(9.2) \quad \tilde{M}_n \otimes \lambda_{[0,1]}^{\otimes \mathbb{N}}((n^{-1/2}T_n^{(k)}, n^{-1/4}\phi_n) \in \cdot) \to \tilde{M}^{(1)} \otimes \lambda_{[0,1]}^{\otimes \mathbb{N}}((T^{(k)}, \phi) \in \cdot)$$

*weakly as probability measures on the space of graph spatial trees. Furthermore,*

$$\lim_{k \to \infty} \limsup_{n \to \infty} \tilde{M}_n \otimes \lambda_{[0,1]}^{\otimes \mathbb{N}}(n^{-1/2}\Delta(\overline{T}_n, \overline{T}_n^{(k)}) > \varepsilon) = 0,$$

*and also*

$$\lim_{\delta \to 0} \limsup_{n \to \infty} \tilde{M}_n \left( \sup_{d_{\overline{T}_n}(\sigma_1, \sigma_2) < \delta n^{1/2}} n^{-1/4} d_E(\phi_n(\sigma_1), \phi_n(\sigma_2)) > \varepsilon \right) = 0,$$

*for every $\varepsilon > 0$.*

(c) *Part* (b) *holds when $\phi_n$, $\phi$ are replaced by $\phi_n^{(k)}$, $\phi^{(k)}$ respectively in* (9.2).

**10. Example: scaling limit for SRW on BRW.** To illustrate the results of the previous sections, we will demonstrate how the simple random walks on the graphs generated by conditioned branching random walks converge to the Brownian motion on the support of the integrated super-Brownian excursion. Let us start by introducing some notation. For an unweighted graph tree $T$ with root $\rho$, let $E^T$ be its edge set, and, for each $\sigma_1, \sigma_2 \in T$, denote by $E^T_{\sigma_1,\sigma_2}$ the subset of $E^T$ containing the $d_T(\sigma_1, \sigma_2)$ edges in the shortest path from $\sigma_1$ to $\sigma_2$ in $T$. Given a function $y : E^T \to \mathbb{R}^d$, we can define a map $\phi : \overline{T} \to \mathbb{R}^d$ by setting $\phi(\rho) := 0$,

$$\phi(\sigma) := \sum_{e \in E^T_{\rho,\sigma}} y(e) \qquad \forall \sigma \in T \setminus \{\rho\}$$

and linearly interpolating along edges. Clearly, $y(e)$ records the increment of $\phi$ along the edge $e \in E^T$.

We can now describe the family of random graph spatial trees

$$\{(T_n, \phi_n)\}_{n \geq 1}$$

that we will consider throughout the remainder of this article. First, for each $n \in \mathbb{N}$, the random ordered graph tree $T_n$ is the family tree generated by a Galton–Watson branching process started from a single ancestor with offspring distribution $Z$ conditioned to have $n$ vertices. Following [17], we will assume that $Z$ satisfies

$$\mathbf{E}Z = 1, \qquad \sigma_Z^2 := \operatorname{Var} Z \in (0, \infty), \qquad \mathbf{E}e^{\lambda Z} < \infty$$

for some $\lambda > 0$. To describe the spatial element of $(T_n, \phi_n)$, we suppose that, conditional on $T_n$, the function $y : E^{T_n} \to \mathbb{R}^d$ is defined so that $(y(e))_{e \in E^{T_n}}$ are independent, each distributed as a random variable $Y$, which is assumed to satisfy

$$\mathbf{E}Y = 0, \qquad \operatorname{Var} Y = \Sigma_Y^2, \qquad \mathbf{P}(d_E(0, Y) \geq x) = o(x^{-4})$$

for some positive definite $d \times d$-matrix $\Sigma_Y$, and then define $\phi_n : \overline{T}_n \to \mathbb{R}^d$ from $y$ as in the previous paragraph. Observe that, conditional on $T_n$, if $\rho, \sigma_1, \ldots, \sigma_l$ is an injective path in $T_n$, then the path $\phi(\rho), \phi(\sigma_1), \ldots, \phi(\sigma_l)$ is a simple random walk in $\mathbb{R}^d$ with step distribution $Y$. Thus, taking into account the independence properties of $y$, it is easy to check that the collection of paths in $\mathbb{R}^d$ obtained by mapping the paths which emanate from the root of $T_n$ into $\mathbb{R}^d$ using $\phi_n$ form a branching random walk, conditioned to have a total of $n$ particles. The one-dimensional version of the following result was proved in [17]; the generalization to $d$ dimensions is straightforward. Note that, given the other assumptions that we are making, [17], Theorem 2, implies that the $o(x^{-4})$ condition on the tail of the distribution of $Y$ is actually necessary to obtain this convergence of tours.

PROPOSITION 10.1 (cf. [17], Theorem 2). *If we define $(v_n, r_n)$ to be the random normalized discrete tour associated with the random graph spatial tree $(T_n, \phi_n)$ for each $n \in \mathbb{N}$, and set*

$$\sigma_T := \frac{2}{\sigma_Z}, \qquad \Sigma_\phi := \Sigma_Y \sqrt{\frac{2}{\sigma_Z}},$$

*then $(v_n, r_n) \to (\sigma_T v, \Sigma_\phi r)$ in distribution in $C([0,1], \mathbb{R}_+) \times C([0,1], \mathbb{R}^d)$, where $(v, r)$ is a random tour with law $\tilde{M}^{(1)}$.*

By rescaling Theorem 9.1 appropriately using $\sigma_T$ and $\Sigma_Y$, we are subsequently able to deduce the convergence of the annealed laws of the simple random walks on $T_n$ mapped into $\mathbb{R}^d$ by $\phi_n$.

COROLLARY 10.2. *If $\tilde{M}_n$ is the law of the random normalized discrete tour $(v_n, r_n)$ associated with $(T_n, \phi_n)$, then $\mathbb{M}_n$, as defined by (9.1), converges to*

$$\mathbb{M}(\{(v, r, \phi(X^{\mathcal{T}})) : (\sigma_T v, \Sigma_\phi r, (\Sigma_\phi \phi(X^{\mathcal{T}}_{t\sigma_T^{-1}}))_{t \geq 0}) \in \cdot\})$$

*weakly as probability measures on $C([0,1], \mathbb{R}_+) \times C([0,1], \mathbb{R}^d) \times C(\mathbb{R}_+, \mathbb{R}^d)$.*

In terms of random variables, this result has the following consequences (recall also Theorem 8.1). We use the notation $\Rightarrow$ to represent convergence in distribution.

COROLLARY 10.3. *If $(v_n, r_n, \phi_n(X^{T_n}))$ has law $\mathbb{M}_n$ and $(v, r, \phi(X^{\mathcal{T}}))$ has law $\mathbb{M}^{(1)}$, then*

$$(n^{-1/2}\overline{T}_n, n^{-1/4}\phi_n) \Rightarrow (\sigma_T \mathcal{T}, \Sigma_\phi \phi),$$

$$n^{-1/4}\phi_n(T_n) \Rightarrow \Sigma_\phi \phi(\mathcal{T}),$$

$$\mu^{S_n}(n^{1/4} \cdot) \Rightarrow \mu^{\mathcal{S}}(\Sigma_\phi^{-1} \cdot),$$

$$(n^{-1/4}\phi_n(X^{T_n}_{tn^{3/2}}))_{t \geq 0} \Rightarrow (\Sigma_\phi \phi(X^{\mathcal{T}}_{t\sigma_T^{-1}}))_{t \geq 0},$$

*simultaneously, in $\mathbb{T}_{\mathrm{sp}}$, as compact subsets of $\mathbb{R}^d$, weakly as Borel probability measures on $\mathbb{R}^d$, and in $C(\mathbb{R}_+, \mathbb{R}^d)$, respectively.*

DEPARTMENT OF STATISTICS
UNIVERSITY OF WARWICK
COVENTRY CV4 7AL
UNITED KINGDOM
E-MAIL: d.a.croydon@warwick.ac.uk